\documentclass{article}
\usepackage{amsfonts}
\usepackage{amssymb}

\begin{document}
\newtheorem{theor}{Theorem}[section] 
\newtheorem{prop}[theor]{Proposition} 
\newtheorem{cor}[theor]{Corollary}
\newtheorem{lemma}[theor]{Lemma}
\newtheorem{sublem}[theor]{Sublemma}
\newtheorem{defin}[theor]{Definition}
\newtheorem{conj}[theor]{Conjecture}

\hfuzz2cm

\def\deg{{\widehat {\rm deg}\,}}
\def\bz{\mbox{\boldmath$\zeta$\unboldmath}}
\def\bzs{\mbox{\boldmath$\zeta'$\unboldmath}}
\def\bzo{\overline{\bz}}
\def\odd{{\rm odd}}
\def\a{\alpha}
\def\bC{{\bf C}}
\def\Td{{\rm Td}}
\def\ch{{\rm ch}}
\def\Hom{{\rm Hom}}
\def\Ad{{\rm Ad}^{1,0}_{G/P}}
\def\KO{H}
\def\K1{K}
\def\covol{{\rm covol}\,}
\gdef\beginProof{\par{\bf Proof: }}
\gdef\endProof{${\bf Q.E.D.}$\par}
\gdef\mtr#1{\overline{#1}}
\gdef\ar#1{\widehat{#1}}

\def \x{\times}
\def \BN{{\bf N}}
\def \BZ{{\bf Z}}
\def \BQ{{\bf Q}}
\def \BR{{\bf R}}
\def \BC{{\bf C}}
\def \BG{{\bf G}}
\def \<{\langle}
\def \>{\rangle}
\def \Spec {\hbox{\rm Spec }}
\def \X {{\rm Bor}(G)}
\def \Y {{\rm Par}(G)}
\def \ad {\hbox {ad}}

\author{ Christian Kaiser\\
Kai K\"ohler} 
\title{A fixed point formula of Lefschetz type in Arakelov geometry
III: representations of Chevalley schemes and heights of flag varieties}
\maketitle
\begin{abstract}
We give a new proof of the Jantzen sum formula for integral representations
of Chevalley schemes over Spec {\bf Z}. This is
done by applying the fixed point formula of Lefschetz type in Arakelov geometry to
generalized flag varieties. Our proof involves the computation of the equivariant
Ray-Singer torsion for all equivariant bundles over complex homogeneous spaces.
Furthermore, we find several explicit formulae for the global height of any generalized
flag variety.
\end{abstract}
\begin{center}
\end{center}
\thispagestyle{empty}
\newpage
\setcounter{page}{1}
\tableofcontents
\newpage
\parindent=0pt
\parskip=5pt

\section{Introduction}
One nice application of the classical Lefschetz fixed point formula is a geometric
interpretation of the Weyl character formula for irreducible representations of a
compact Lie group $G_c$, found by Atiyah and Bott \cite{AB}. Namely, choose a
maximal torus
$T_c$ of
$G_c$, fix a weight
$\lambda$ and consider the line bundle associated to the $T_c$-representation of
weight
$\lambda$ over the flag variety $G_c/T_c$. The action of $G_c$ on $G_c/T_c$ induces
an action on the cohomology vector spaces of this line bundle. The Lefschetz fixed
point formula provides a formula for the character of this $G_c$-representation as a
sum over the fixed point set. The fixed point set can be identified with the Weyl
group of $G_c$, and the resulting formula is the classical Weyl character formula.

The purpose of this article is to investigate the analogous problem in the context of
Arakelov geometry. Instead of compact Lie groups, we consider Chevalley schemes $G$
over Spec {\bf Z}. For a technical reason, we have to exclude factors of type $G_2,
F_4$ and
$E_8$. The generalized flag varieties associated to $G$ are acted upon by a maximal
split torus of
$G$. We apply to them the analogue of the fixed point formula in Arakelov geometry
\cite{KR1} due to Roessler and one of the authors. If
$V_\mu$ denotes the weight space to the weight
$\mu$ in the cohomology representation, then the classical Lefschetz formula computes
the character
$
\sum_{\mu\ {\rm weight}} \mu \dim V_\mu\,\,.
$
In the Arakelov situation, the cohomology is a $\bf Z$-module whose tensor product
with $\bf C$ carries an Hermitian metric, and the fixed point
formula computes the arithmetic character
$$
\sum_{\mu\ {\rm weight}} \mu \log \frac{{\rm covol}\,V_{{\rm free},\mu}}{\#\,V_{{\rm
tor},\mu}}\,\,.
$$
with $V_{{\rm free}},\,V_{\rm
tor}$ indicate the free and torsion part, respectively.
Having chosen a maximal split torus $T\subseteq G$ we have a canonical maximal
compact subgroup $G_c$ in $G({\bf C})$. Our metric on the cohomology is 
$G_c$-invariant and thus determined up to a factor. The arithmetic character of the
alternating sum of the cohomology modules is then expressed in terms of roots and the
Weyl group of $G_c$.

The resulting formula (up to the determination of the $L^2$-metric) turns out to be the
Jantzen sum formula
\cite[p. 311]{Jan} which has been proven originally using solely methods from
representation theory and algebraic geometry. In contrast, in the proof presented here
we mostly have to deal with the differential geometry of Laplace operators and we use
only the representation theory of compact Lie groups. The arithmetic geometry is reduced
to that of the fixed point scheme which consists of copies of Spec $\bf Z$.

We have to compute the equivariant
holomorphic Ray-Singer torsion associated to vector bundles on complex homogeneous
spaces. This is based on the formula for the Ray-Singer torsion for Hermitian
symmetric spaces which has been determined by explicitly describing the zeta function
defining the torsion by one of the authors \cite{K2}. Using the
arithmetic Lefschetz formula or a formula due to Ma \cite{Ma2}, we can compare the
analytic torsion of the base, of the total space and of the fiber of a fibration. Now we
consider a tower of fibrations of complex homogeneous spaces whose fibers are Hermitian
symmetric. The lowest base is a point and the total space is the full flag space. By an
induction procedure, we get a formula for the torsion on any full flag space. Using
again one fibration over an arbitrary complex homogeneous base and with total space and
fiber being full flag spaces, we get an expression for the analytic torsion for that
base.

In the last chapter, we use the Jantzen sum formula for the arithmetic character to
derive formulae for the global height $h(X,\mtr L)$ of ample line bundles $\mtr L$ over
any (generalized) flag variety
$X$ (now including the types $G_2$, $F_4$, $E_8$). The global height is defined using
characteristic classes in Arakelov geometry. 

The global height of generalized flag varieties has already been investigated in
numerous cases using Arakelov intersection theory. There are explicit formulae for
projective spaces due to Bost, Gillet and Soul\'e \cite{BoGS} and for even-dimensional
quadrics by Cassaigne and Maillot \cite[Cor. 2.2.10]{CM}. Maillot \cite{M} and
Tamvakis
\cite{T1} found algorithms for the height of Grassmannians, leading to an explicit
formula in the case
$G(2,n)$ \cite{T1}. Tamvakis also found combinatorial algorithms giving the height of
generalized flag varieties of type $A_n$ \cite{T2} and of Lagrangian
Grassmannians
\cite{T3}.

Let Ht denote the additive topological characteristic class associated to the power
series
$$
{\rm Ht}(x):=\sum_{k=0}^\infty \frac{(-x)^k}{2(k+1)(k+1)!}\,\,.
$$
This is the Taylor expansion of the function
$x\mapsto\frac1{2x}\left(\log|x|-\Gamma'(1)-{\rm Ei}(-x)\right)$ at $x=0$ with Ei
being the exponential-integral function. Let $\mtr {\cal L}$ be an equivariant line
bundle on
$G/P$, ample over $\bf C$.
We construct in a simple
way a virtual holomorphic vector bundle
$y$ on
$G_{\bf C}/P_{\bf C}$ such that the height with respect to 
$\mtr {\cal L}$ is given by
$$
h(G/P,\mtr {\cal L})=(n+1)!\int_{G_{\bf C}/P_{\bf C}}
{\rm Ht}(y) e^{c_1(L)}
$$
(Theorem \ref{height}). Also we give explicit formulae for the height in terms of roots
and the Weyl group of
$G_c$. From these formulae, we derive some simple properties of heights on
generalized flag varieties.

Our method to compute the height is similar to the approach given in \cite{KR2} by
Roessler and one of the authors. Still, it is only an application of the classical
Jantzen sum formula combined with a relation between the height and the asymptotic of
the covolume of the cohomology of
$\mtr L^m$ for large
$m\in{\bf N}$ which is described in
\cite[VIII.2.3]{Soule}. The proof does not make direct use of the Lefschetz
fixed point formula in Arakelov geometry.

If one considers just the special case of Hermitian symmetric spaces, chapters 3-6 are
not necessary for the proofs in chapters 7 and 8. In particular the computation
of the analytic torsion for generalized flag varieties in these chapters is already
known in the case of Hermitian symmetric spaces (\cite{K2}). 

Furthermore, the full main result of \cite{K2} is not needed in this article. In
chapters 3-6 we apply it only in the very special case of actions with isolated fixed
points. In this case, the proof in
\cite{K2} might be replaced by a considerably shorter one similar to the
proof of \cite[Lemma 13]{K3}.

 This article is a part of the habilitation
thesis of the second author.

{\bf Acknowledgements} The authors wish to thank Gerd Faltings, G\"unter Harder,
Jens-Carsten Jantzen, Damian
Roessler and Matthias Weber for many valuable discussions and comments. We are also very
grateful to the referees for their detailed comments.

\section{Flag Varieties}

In this chapter we recall the definition and some properties of the
flag variety of a Chevalley group over $\BZ$. Further on we describe
the fixed point scheme under the action of some multiplicative subgroup
scheme.

Notations: We denote the category of schemes resp. sets by
$\underline{\hbox{Scheme}}$ resp. $\underline{\hbox{Set}}$. For any scheme
$S$ we denote by
$\underline{\hbox{Scheme}_S}$
the category of $S$-schemes.
We identify a $S$-scheme by the functor from
$\underline{\hbox{Scheme}_S}$  to $\underline{\hbox{Set}}$ which it represents.
If $X$ and $S$ are schemes we set $X_S:=X\x S$ and consider it as an
$S$-scheme.  If $S=\Spec A$ is affine we also write $X_A$. For $A=\BC$ we also denote by $X_\BC$
the complex analytic variety $X(\BC)$. For a $\BZ$-module
$M$ and ring $A$ we set $M_A:=M\otimes A$. Lie algebras associated to smooth group schemes shall
be denoted by the corresponding german letter.
\par Let $G$ be a semisimple Chevalley group of rank $r$.
This is a smooth affine group scheme over
$\Spec \BZ$ with connected semisimple groups as geometric fibres. 
A closed subgroup
scheme of $G_S$ is called a Borel subgroup, resp. a parabolic subgroup, if
it is smooth and of finite presentation over $S$ and a Borel subgroup,
resp. a parabolic subgroup,
in every geometric fibre.

Consider the functors from $\underline{\hbox{Scheme}}$ to
$\underline{\hbox{Set}}$ which map $S$ to the set of all Borel subgroups, resp.
parabolic subgroups, of $G_S$.
These functors are representable by smooth projective schemes $\X$, resp.
$\Y$, over $\Spec \BZ$ \cite[ch. XXVI, 3.3]{SGA3}.
The group $G$ acts by conjugation on these schemes.

We fix now a maximal split torus $T\subseteq G$ with group of characters $X^*(T)$
and cocharacters
$X_*(T)$, set of roots
$\Sigma$ and Weyl group (of the root system) $W_G$. The adjoint action of $T$
on
$\frak g$ gives the weight decomposition in free $\BZ$-modules:
$$\frak g =\frak t \oplus \bigoplus_{\alpha\in \Sigma} \frak g_{\alpha}.
$$
Consider a parabolic subgroup $P$ containing $T$. Since the root spaces $\frak
g_{\alpha}$ are of rank one there is a subset $R\subseteq \Sigma$ with
$$\frak p =\frak t \oplus \bigoplus_{\alpha\in R} \frak g_{\alpha}.
$$
The subset $R\subset\Sigma$ is closed (i.e. if $\alpha, \beta\in R$ and
$\alpha +\beta\in
\Sigma$ then $\alpha +\beta\in  R$) and
$R\cup
-R=\Sigma$ . We call such a subset of roots a
parabolic subset.
\medskip
\begin{lemma}\label{ch11}
 For each parabolic subset of roots
$R\subset \Sigma$ there exists a unique parabolic subgroup $P(R)$
in
$G$ containing $T$ such that $\frak p(R)=
\frak t \oplus \bigoplus_{\alpha\in R} \frak g_{\alpha}$.
\end{lemma}
\beginProof The same proof as in \cite[XXII, 5.5.1]{SGA3}. \endProof
\medskip
For a $S$-group scheme $H$ acting on a $S$-scheme $Y$ we denote by $Y^H$ the
functor of fixed points, i.e. for a $S$-scheme $S^\prime$ we have
$$Y^H(S^{\prime}):=\{y\in Y(S^{\prime})\;\vert\;\forall \phi
:S^{\prime\prime}\to S^{\prime}, \;\forall
h\in H(S^{\prime\prime}):\; h\phi^* (y)=\phi^*(y)\}.
$$
In the braces $\phi$ is an arbitrary
$S$-morphism and $\phi^*$ the induced map on points.
 \medskip
\begin{prop}\label{ch12}
The functor morphism $\Y^T\to \Y$ is
represented by the closed embedding of copies of $\Spec \BZ$ given by
$P(R)\in \Y(\BZ)$:
$$ \bigsqcup_{R\subseteq \Sigma\hbox{ \rm parabolic} } \Spec \BZ\to \Y.
$$
\end{prop}
\beginProof Let $P$ be a parabolic subgroup over $S$. Then $P$
is its own normalizer (\cite[XXII, 5.8.5]{SGA3}),
i.e. for any $S$-scheme $S^{\prime}$ and any $g\in G(S^{\prime})$ with
$\hbox{int}(g)(P_{S'})=P_{S'}$ we have
$g\in P(S^{\prime})$. Hence we have $P\in
\Y^T(S)$ iff
$T_S\subset P$.
The parabolic subgroups containing $T_S$ are determined by their Lie
algebra
(\cite[XXII,
5.3.5]{SGA3}). For $S$ connected we have $\frak p=
\frak t\otimes {\cal O}_S \oplus \bigoplus_{\alpha\in R} \frak g_{\alpha}\otimes {\cal O}_S$ for some parabolic
subset of roots (\cite[XXVI, 1.4]{SGA3}), and hence $P=P(R)_S$.
Since clearly $\Y^T(S_1\sqcup S_2)= \Y^T(S_1)\x \Y^T(S_2)$ the claim follows.
\endProof
\medskip
Let us fix an ordering $\Sigma=\Sigma^+\cup\Sigma^-$ with base $\Pi$.
A parabolic subset $R$ containing $\Sigma^-$, resp. its associated parabolic
subgroup $P(R)$, is called "standard". For a standard parabolic subset $R$ we
set
$\hbox{type}(R)= \{\alpha\in \Pi\;\vert\;\alpha\in R\cap -R\}$.
For $k$ a field every parabolic subgroup $P$ of $G_k$ is conjugate to an
unique
standard parabolic $P(R)$. We set $\hbox{type}(P):= \hbox{type}(R)$.
%
For general $S$ it is a locally
constant function
$s\to \hbox{type}(P_s)$.

We have a type
decomposition in open and closed subsets:
$$\Y=\bigsqcup_{\Theta\subseteq \Pi}\Y_{\Theta}
$$
with $\Y_{\Theta}$ classifying parabolic subgroups of type $\Theta$.
Two parabolic groups have the same type if and only
if they are locally conjugate w.r.t. the fpqc or \'etale topology (\cite[XXVI,
3.3]{SGA3}). Since
parabolic subgroups are their own normalizer we see that
$\Y_{\Theta}$ represent also the fpqc- or \'etale-sheafification of the
functor
$S\to G(S)/P(S)$ for any $P\subseteq G$ of type $\Theta$. We write therefore
$\Y_{\Theta}=:G/P$. It is a connected scheme. We have
$\Y_{\emptyset}=\X$ and
$\Y_{\Pi}=\Spec \BZ$. 
Denoting the standard parabolic of type $\Theta$ by $P_\Theta$, 
Proposition \ref{ch12} implies
 
\begin{cor}\label{ch13}
 For $\Theta\subseteq\Pi$
we have
$$(G/P_{\Theta})^T=
\bigsqcup_{w\in W_G/W_{\Theta} } \Spec \BZ.
$$
Here $W_{\Theta}$ is the subgroup generated by the
reflections at roots in $\Theta$.
\medskip
\end{cor}
\par\noindent {\bf Remark }: By \cite[XXVI, 3.16 (iii)]{SGA3} the fixed point scheme
$(G/P)^T$ is finite \'etale over $\Spec \BZ$, from that Corollary \ref{ch13}
also follows.

We consider now a finite diagonalisable closed subgroup scheme $\mu\subseteq T$.
This means we fix a finite abelian group $M$ and a surjective map
$\psi: X^*(T)\to M$. Then $\mu\simeq\Spec \BZ[M]$ with comultiplication
given by $\BZ[M]\to\BZ[M]\otimes \BZ[M],\; m\mapsto
m\otimes m$. We assume now the following property for $\mu$:
\medskip
\begin{defin} \label{propReg}
 A
finite diagonalisable closed subgroup scheme $\mu\subseteq T$ has the property
$(\hbox{Reg})$ iff the restriction of $\psi$ to the set of roots $\Sigma$ is injective.
\end{defin}
\begin{prop} If $\mu\subseteq T$ has property
$(\hbox{Reg})$ then $\Y^{\mu}= \Y^{T}$.
\end{prop}
\beginProof Take $P\in \Y(S)$. We have to show that $P$
is fixed by $\mu_S$ iff it is fixed by $T_S$. For this we can assume $G$ to be
adjoint.
From \cite[XXII, 5.3.3]{SGA3} it follows that $P$ is invariant under the action
of $T_S$ resp. $\mu_S$ iff the same is true for $\frak p$.
It is a general property of a diagonalisable
group scheme
acting linearly on a module $M$ with weight decomposition
$M=\bigoplus_{\alpha} M_{\alpha}$, that a submodule
$N\subseteq M$ is
invariant iff $N=\bigoplus_{\alpha} (N\cap M_{\alpha})$. Since by
property $(\hbox{Reg})$ the
weight decompositions of $\frak g$ with respect to $\mu$ and $T$
coincide the claim follows. \endProof
\medskip
Let $\<\cdot,\cdot\>:X^*(T)\x X_*(T)\to\BZ$ denote the natural pairing. For a coroot
$\alpha^{\lor}$ we set $H_{\alpha}:=d\alpha^{\lor}(1)\in \frak t$.
For each $\alpha\in\Sigma$ there is a homomorphism $x_{\alpha}:\BG_a\to G$ such
that $dx_{\alpha}:{\rm Lie\,}\BG_a\to \frak g_{\alpha}$ induces an isomorphism. This
map is unique up to a sign change. We can choose the maps $x_\a$ in such a way that
with
$X_{\alpha}:=dx_{\alpha}(1)$ we get a Chevalley base $\{X_{\alpha},
H_{\mu}\;\vert\;\alpha\in \Sigma,\mu\in \Pi\}$ of $\frak g\otimes\BC$. For $G$
simply connected it is actually a base of $\frak g$.
There is a unique involution $\tau$ acting on $G$ (\cite[Corollary II.1.16 and
proof]{Jan}) with:
\par 1) ${\tau}_{|T}: t\mapsto t^{-1}$
\hskip 20pt 2) $\tau\circ x_{\alpha}(\cdot)=x_{-\alpha}(-\cdot)$.
\par\noindent
It is independent from the consistent choice above.

We can now define a maximal compact subgroup $G_c$ of $G_\BC$:
$$ G_c:=(G_\BC)^{\tau\circ \iota}.
$$
Here $\iota$ denotes the complex conjugation.
\medskip
\begin{lemma}\label{ch16}
 The group $N_G(T)(\BZ)\cap G_c$ acts on
$T$ by conjugation through the full Weyl group.
\end{lemma}
\beginProof For each root $\alpha$ set
$n_{\alpha}=
x_{\alpha}(1) x_{-\alpha}(-1) x_{\alpha}(1)$. By \cite[p.176]{Jan} we have
$n_{\alpha}\in
N_G(T)(\BZ)$ and it acts on $T$ by the reflection at $\alpha$. Since $\tau
(x_{\alpha}(1) x_{-\alpha}(-1)
x_{\alpha}(1))= x_{-\alpha}(-1) x_{\alpha}(1) x_{-\alpha}(-1)$ it is only an
easy calculation in $SL_2$ to check
$\tau(n_{\alpha})=n_{\alpha}$ \cite[p.176]{Jan}.
\endProof
\medskip
We can use the maximal compact group $G_c$ to endow the complex flag variety
$(G/P)_\BC$ with a ca\-no\-ni\-cal hermitian metric which fits nicely with
the integral structure. Fix a standard parabolic subgroup $P\subseteq G$. 
Setting $K:=P_\BC\cap G_c$ the inclusion induces
an isomorphism $(\frak g/\frak p)\otimes \BC\simeq \frak g_c/\frak k$. For $\a$ a
weight in $\frak g/\frak p$ it maps the integral generator of the $\a$-weight space
$X_{\alpha}$ to
$X_{\alpha}- X_{-\alpha}$.
The negative of the Killing form $\kappa(\cdot,\cdot)$ induces
a positive definite real bilinear form $(\cdot,\cdot)$ on $\frak m:=\frak g_c/\frak
k$. Denote the complex structure on $\frak m$ by $J$. Set
$\frak m^{1,0}:=\{X\in\frak m\otimes \BC\,|\,J X=i X\}$ and define $X^{1,0}:=(X-i J
X)/2$ for
$X\in\frak m$. The metric on $\frak m$ induces a $G_c$-invariant Hermitian metric on
$\frak m^{1,0}$ and thus an Hermitian metric $(\cdot,\cdot)_B$ on $(G/P)_\BC$. Using
the Killing form we get a norm on 
$(\frak t\otimes \BR)^\vee\cong X^*(T)\otimes \BR$. We get (using
\cite[p.37]{Hu})
$$(X_{\alpha},X_{\alpha})=-\kappa(
X_{\alpha}- X_{-\alpha}, X_{\alpha}- X_{-\alpha})
=2\kappa (X_{\alpha},X_{-\alpha})={4\over \|\a\|^2}
$$
and thus $\|X_\alpha^{1,0}\|_B^2=\frac2{\|\a\|^2}$.
By corollary \ref{ch13} the normal bundle of $(G/P)^T$ is canonically
isomorphic to the restriction of the relative tangent bundle $Tf$ to $(G/P)^T$,
where		
$f:G/P\to\Spec\BZ$ is the structure morphism. Using Lemma \ref{ch16} 
and denoting by
$\Psi\subseteq
\Sigma^+$ the set of weights in
$\frak g/\frak p$, the above calculation gives
\medskip
\begin{lemma} \label{ch17}
For $\a\in\Psi$ the generators of the $\alpha$-weight spaces
of the normal bundle of $(G/P)^T$ have squarelength
${2\over \|\a\|^2}$ on every component.
\end{lemma}
Let $P$ be a standard parabolic subgroup.
For $\lambda\in X^*(P)$ we define an invertible sheaf ${\cal L}_{\lambda}$
on
$G/P$ by
$${\cal L}_{\lambda}(U)=\bigl\lbrace f\in { \cal O}( \pi^{-1}U)\;\vert\;
  f(xp)=\lambda(p)^{-1}f(x)\;\; \forall x\in (\pi^{-1}U)(S),p\in P(S)\hbox{
and
}\forall S \bigl\rbrace.
$$
Here $\pi:G\to G/P$ is the canonical map and $U\subseteq G/P$ an open
subset.
Considering $X^*(P)$ as a
subgroup of $X^*(B)$ we have a canonical isomorphism
$H^*(G/P,{\cal L}_{\lambda})\simeq H^*(G/B,{\cal L}_{\lambda})$. Hence
in the following we consider only the "full flag variety". The groups
$H^*(\lambda):=H^*(G/B,{\cal L}_{\lambda})$ are finitely generated
$\BZ$-modules with an algebraic action of $G$.
We use the "dotted
action" of the Weyl group: $w\ldotp
\mu=w(\mu+\rho)-\rho$. Here $\displaystyle \rho={1\over
2}\sum_{\alpha\in \Sigma^+}\alpha$. Let $I_+$ denote the set of dominant weights
w.r.t. $\Sigma^+$. For the action of $G_{\BQ}$ on
$H^*(w\ldotp\lambda)_{\BQ}$ we have, denoting the
irreducible representation of highest weight $\lambda$ by $V_{\rho+\lambda}$ and by
$l(w)$ the length of a $w\in W_G$:
\medskip
\begin{theor} \label{ch18}
 (Borel-Weil-Bott): a) If $\lambda+\rho\in I_+$ 
but
$\lambda\notin I_+$, then
$$H^*(w\ldotp\lambda)_{\BQ}=0\hskip1.5cm \hbox{for all }w\in W_G.
$$
\par \noindent b)
If $\lambda\in I_+$ then we have for all $w\in W_G$ and all $i\in
\BN_0$:
$$ H^i(w\ldotp\lambda)_{\BQ}\simeq\left\{
\begin{array}{cc}
V_{\rho+\lambda}&\hskip1.5cm\hbox{if }i=l(w)\cr
0\; &\hskip1.6cm \hbox{otherwise.}
\end{array}\right.
$$
\end{theor}
\bigskip
For an abelian group $M$ we denote by $M_{\rm tor}$ its torsion and set
$M_{\rm free}=M/M_{\rm tor}$. Then $H^i(w\ldotp\lambda)_{\rm free}$ is a lattice in
$H^i(w\ldotp\lambda)_{\BQ}$.
The free part of the following Proposition is due to Andersen [A1, 2.10], the
torsion part is a remark of Faltings. We set $X=G/B$ and $n+1=\dim X$.

\begin{prop} \label{Serre}
Let $\lambda$ be an arbitrary weight.
For all $i\in{\bf Z}$ there are perfect pairings of $G$-modules
$$ H^i(\lambda)_{\rm free}\x H^{n-i}(-\lambda -2\rho)_{\rm free}\to {\bf Z} $$
$$ H^i(\lambda)_{\rm tor}\x H^{n+1-i}(-\lambda -2\rho)_{\rm tor}\to {\bf Q/Z}. $$
\end{prop}
\beginProof Applying Grothendieck duality \cite{Hart2} to the smooth and
proper structure morphism $f:X\to \Spec \bf Z$ gives an isomorphism in the
derived category of bounded complexes of abelian groups
$$ R\Gamma ({\cal L}_{-\lambda-2\rho})[n]\stackrel\sim\to R\Hom_{\bf Z}
(R\Gamma ({\cal L}_{\lambda}),\bf Z).
$$
The right hand side is a composition of derived functors and hence the
limit
of a (in this case: second quadrant) spectral sequence $(E_r^{i,j})$, with
$E_2^{i,j}=\hbox{Ext}^j_{\bf Z}(H^{-i}(\lambda),\bf Z)$. Since
$\hbox{Ext}^j_{\bf Z}(\cdot,\cdot)=0 \;\forall j\geq 2$
the spectral sequence degenerates at $E_2$, and we get exact sequences
$$0\to \hbox{Ext}^1_{\bf Z}(H^j(\lambda),{\bf Z})\to
H^{n+1-j}(-\lambda-2\rho) \to \Hom (H^{j-1}(\lambda),\bf Z)\to 0 .
$$
The group $\hbox{Ext}^1_{\bf Z}(H^j(\lambda),{\bf Z})=\hbox{Ext}^1_{\bf Z}
(H^j(\lambda)_{\rm tor},{\bf Z})\simeq \Hom(H^j(\lambda)_{\rm tor},{\bf Q/Z})$ is
torsion and $\Hom (H^{j-1}(\lambda),\bf Z)$ is torsionfree. The claim follows.
\endProof

\medskip
For a $T$-module $A$, which is a finite abelian group, we set $\hbox{\rm char
}A:=-\sum_{\mu}\mu \log \#A_{\mu}$.

\begin{cor} \label{Serretors}
$\hbox{\rm char }H^i(\lambda)_{\rm tor}= \hbox{\rm char
}H^{n+1-i}(w_0.\lambda)_{\rm tor}$.
\end{cor} 
\beginProof
For $n_0\in N_G(T)(\BZ)$ acting on $T$ by $w_0$
we consider the automorphism
$\widetilde{\tau}:=
\hbox{Int}(n_0)\circ \tau$. This automorphism respects
$B$, acts therefore on $G/B$, and puts ${\cal L}_{\lambda}$ to ${\cal
L}_{-w_0\lambda}$. Hence it induces an isomorphism
$H^i(\lambda)\stackrel\sim\to H^i(-w_0\lambda)$, giving the equality
$\hbox{\rm char }\Hom(H^i(\lambda)_{\rm tor},{\bf Q/Z}) =\hbox{\rm char
}H^i(-w_0\lambda)_{\rm tor}$. The claim follows from the Proposition.
\endProof

We assume in the following $\lambda$ to be dominant.
 We choose a generator $v_{w\ldotp\lambda}$ of the
$\lambda$-weight space of $H^{l(w)}(w\ldotp\lambda)_{\rm free}$. By
irreducibility there is a unique $G_c$-invariant hermitian metric
$(\cdot,\cdot)_{\tau}$ on
$H^{l(w)}(w\ldotp\lambda)_{\BC}$
with
$(v_{w\ldotp\lambda},v_{w\ldotp\lambda})_{\tau}=1$. Its restriction to
$H^{l(w)}(w\ldotp\lambda)_{\BR}$ is a $\tau$-contravariant form (\cite[p. 35]{J1}),
i.e. $(gv,v^{\prime})_{\tau}= (v,\tau(g)^{-1}v^{\prime})_{\tau}$,
and was studied by Jantzen \cite{J1} and Andersen \cite[p. 515]{A1}. By
theorem 2.8 above we have for any $w\in W$ and for $w_0\in W$ the longest element
an unique isomorphism:
$$T_{w_0}:H^{l(w)}(w\ldotp\lambda)_{\BQ}\to
H^{l(w_0w)}(w_0w\ldotp\lambda)_{\BQ}
$$
with $T_{w_0}(v_{w\ldotp\lambda})= v_{w_0w\ldotp\lambda}$.
\medskip
\begin{prop} \label{ch19}
The lattices
$T_{w_0}(H^{l(w)}(w\ldotp\lambda)_{\rm free})$ and $ H^{l(w_0w)}(w_0w\ldotp\lambda)_{\rm free}$
are dual w.r.t. the pairing
$(\cdot,\cdot)_{\tau}$.
\end{prop}
\beginProof This follows immediately from the explicit
presentation of a contravariant form in \cite[p.515]{A1}. For convenience of
the reader we recall it. 
For $n_0\in N_G(T)(\BZ)$ acting on $T$ by $w_0$
consider as above
$\widetilde{\tau}:=
\hbox{Int}(n_0)\circ \tau$, inducing an isomorphism
$\phi:H^{l(w)}(w\ldotp\lambda)\stackrel\sim\to
H^{l(w)}(-w_0(w\ldotp\lambda))$. We set
$\psi:=n_0^{-1}\circ \phi$. It is independent of the choice of $n_0$ and
$\psi(gv)= \tau(g)\psi(v)$. Using Serre duality pairing,  for $v,v^{\prime}\in
H^{l(w)}(w\ldotp\lambda)_\BQ$ we set :
$$ \beta (v,v^{\prime})=(T_{w_0}(v),\psi(v^{\prime})).$$
Note that $\beta(v_{w\ldotp\lambda}, v_{w\ldotp\lambda})=\pm 1$.
We have
$\beta(gv,v^{\prime})= \beta(v,\tau(g)^{-1}v^{\prime})$.
Hence it defines a non trivial $\tau$-con\-tra\-va\-ri\-ant form on
$H^{l(w)}(w\ldotp\lambda)_{\BQ}$. We denote by $T_{w_0*}(\beta)$
the transported form on $H^{l(w_0w)}(w_0w\ldotp\lambda)_{\BQ}$. Then
w.r.t. this form
$H^{l(w_0w)}(w_0w\ldotp\lambda)_{\rm free}$ and $T_{w_0}
(H^{l(w)}(w\ldotp\lambda)_{\rm free})$ are dual lattices, as Serre duality induces a
perfect pairing on $H^\bullet_{\rm free}$ (Proposition \ref{Serre}). On the other hand
by uniqueness of normalized $\tau$-contravariant forms we have $T_{w_0*}(\beta)=
\pm (\cdot,\cdot)_{\tau}$. \endProof
\medskip
Take $w_{\lambda}$ in the $\lambda$-weight space of $V_{\rho+\lambda}$ and
$w_{-\lambda}$ in the $-\lambda$-weight space of $V_{\rho+\lambda}^\vee$ with
$\<w_{\lambda}, w_{-\lambda}\>=1$, and consider the associated matrix
coefficient
$c_{w_{\lambda}, w_{-\lambda}}(g)=\<w_{\lambda}, gw_{-\lambda}\>$. The following result
shall be used to refine our results in the case of positive line bundles:
\begin{prop} \label{ch110}
The function
$c_{w_{\lambda}, w_{-\lambda}}$ is a generator of the $\lambda$-weight space
of $H^0(\lambda)$.
\end{prop}
\beginProof It is clear that $c_{w_{\lambda}, w_{-\lambda}}$
lies in the $\lambda$-weight space of $H^0(\lambda)_{\BQ}$.
Since $c_{w_{\lambda}, w_{-\lambda}}(1)=1$ it is enough to show that it
defines a regular function on $G$. We can do this in two ways:
\par 1) We can realize $V_{\rho+\lambda},V_{\rho+\lambda}^\vee$ and the natural
pairing between them by $H^0(\lambda)_{\BQ}$, $H^{l(w_0)}(-\lambda-2\rho)_{\BQ}$ and
the Serre duality pairing. By Proposition \ref{Serre} this defines already a perfect
pairing over $\BZ$. Hence we can take appropriate elements
in the integral structure for
$w_{\lambda}$ and $ w_{-\lambda}$. Therefore $c_{w_{\lambda}, w_{-\lambda}}$
is a regular function on $G$.
\par 2) Using \cite[Proposition I.10.12]{Jan} we have to show that
$\mu( c_{w_{\lambda}, w_{-\lambda}})\in \BZ$ for all distributions $\mu\in
\hbox{Dist}(G)$. Almost
by definition we have $\mu(
c_{w_{\lambda}, w_{-\lambda}})=\<\mu w_{\lambda}, w_{-\lambda}\>$.
By \cite[p. 185]{Jan}  $\hbox{Dist}(G)$ is the subalgebra of
$U(\frak g)$ generated by $X_{\alpha}^n/n!$ with $\alpha\in \Sigma$ and $n\in \BN$, and
by all ${H\choose m}$ with $H\in \frak t$ and $m\in \BN$. There is an
integral PBW-decomposition $U(\frak g)=U^-U^0U^+$ with
$U^-=\<X_{\alpha}^n/n!\;\vert\; \alpha\in \Sigma^-,\;n\in\BN\>$,
$U^+=\<X_{\alpha}^n/n!\;\vert\; \alpha\in \Sigma^+,\;n\in\BN\>$ and
$U^0=\< {H\choose m} \;\vert\;H\in \frak t,m\in \BN\>$.
Since the
vectors
$w_{\lambda}$ and $ w_{-\lambda}$ are highest weight vectors
w.r.t. $B$ resp. the opposite Borel subgroup $B^-$ one has only to consider
$\mu\in U^0$. Considering
$U^0$ as a polynomial algebra on $X^*(T)$ we get
$\mu( c_{w_{\lambda}, w_{-\lambda}})=\mu(\lambda)$.
Since $\frak t=\{ d\phi(1)\;\vert\;\phi\in X_*(T)\}$ we see that the
polynomials ${H\choose m}$ take integral values on $X_*(T)$. The claim
follows.
\endProof

\section{Equivariant Ray-Singer torsion}

In this section we repeat the definition of the equivariant Ray-Singer analytic
torsion introduced in \cite{K1}. Let
$M$ be a compact
$n$-dimensional hermitian manifold with associated K\"ahler form $\omega$. Let $\mtr E$
denote an hermitian holomorphic vector bundle
on $M$ and let
$$
\bar{\partial}:\Gamma^\infty(\Lambda^qT^{*0,1}M\otimes E)
\to\Gamma^\infty(\Lambda^{q+1}T^{*0,1}M\otimes E)
$$
be the Dolbeault operator. As in
\cite{GSZ}, we equip
${\frak A}^{0,q}(M,E):=\Gamma^\infty(\Lambda^qT^{*0,1}M\otimes E)$ with the hermitian 
$L^2$-metric
\begin{equation}\label{metric}
(\eta,\eta'):=\int_M\langle\eta(x),\eta'(x)\rangle\frac{\omega^{\wedge
n}}{(2\pi)^n n!}.
\end{equation} Let $\bar{\partial}^*$ be the adjoint of $\bar{\partial}$
relative to this metric and let
$\square_q:=(\bar{\partial}+\bar{\partial}^*)^2$ be the
Kodaira-Laplacian acting on $\Gamma^\infty(\Lambda^qT^{*0,1}M\otimes E)$ with spectrum
$\sigma(\square_q)$. We denote by ${\rm Eig}_\lambda(\square_q)$ the eigenspace of $\square_q$
corresponding to an eigenvalue $\lambda$. Consider a holomorphic
isometry $g$ of $M$ and assume given a holomorphic isometry $g^E:\mtr E\to g^*\mtr E$.
The fixed point set of $g$ shall be denoted by $M^g$. This notation is chosen different
from \cite{KR1}, \cite{KR2} to avoid other notational problems. The Dolbeault cohomology
$H^{0,q}(M,E):=\ker
\square_q$ shall be equipped with the restriction of the $L^2$-metric. The element
$g$ induces an isometry
$g^*$ of 
$H^{0,q}(M,E)$. Then the equivariant Quillen
metric is defined via the zeta function
$$
Z_g(s):=\sum_{q>0}(-1)^{q+1} q\sum_{\lambda\in\sigma
(\square_q) \atop \lambda\neq 0}
\lambda^{-s}\,{\rm Tr }\,g^*_ {|{\rm Eig}_\lambda(\square_q)}
$$
for ${\rm Re}\,s\gg 0$. Classically, this zeta function
has a meromorphic continuation to the complex plane which is holomorphic
at zero (\cite{Donn}).
\begin{defin} Set $\lambda_g(M,E):={\det}_g
H^{0,*}(M,E):=\bigotimes_q ({\det}_g
H^{0,2 q}(M,E)\otimes{\det}_g
H^{0,2 q+1}(M,E)^\vee)$. The equivariant analytic torsion \cite{K1} is defined as
$$
T_g(M,\mtr E):=Z'_g(0)\,\,.
$$
\end{defin}
For group actions with isolated fixed points on $M$, we consider the equivariant
characteristic classes defined on $M^g$ which equal on every fixed point $p$
\begin{equation}\label{TdDef}
\Td_g(E)_{|p}=\det(1-(g^E_{|p})^{-1})^{-1}
\end{equation}
and
$$
\ch_g(E)_{|p}={\rm Tr }\, g^E_{|p}\,\,.
$$
These classes are defined in a more general context in \cite[Th. 3.5]{KR1}.

\section{A submersion formula for the analytic torsion}\label{masection}
In this section we state a special case of a result due to Xiaonan Ma \cite{Ma2},
namely for isolated fixed points. In this special case, it is also a direct
consequence of the arithmetic Lefschetz formula in its formulation in \cite[end of
section 2, (b)]{KR2}. We refer to
\cite[sections 3, 4 and 6.2]{KR1} for notation and definitions.

Consider a holomorphic K\"ahler fibration $f:M\to B$ of K\"ahler manifolds with fibre
$Z$. Assume given an isometric automorphism $g$ of $f$ such that
all fixed points are isolated. Take an equivariant holomorphic Hermitian vector bundle
$\mtr E$ on $M$ and assume that the dimension of the cohomology groups
$H^*(Z,E_{|Z})$ does not vary on $B$. For
$\zeta\in S^1$ the short exact sequence
$$
0\to TZ\to TM\to f^*TB\to0
$$
decomposes over the fixed point set $M^g$ into $\zeta$-eigenspaces with respect
to the action of
$g$
$$
0\to TZ_\zeta\to TM_\zeta\to f^*TB_\zeta\to0\,\,.
$$
Let $s_\zeta$ denote a non-zero element in $\det TM_\zeta$ and let $s'_\zeta\in\det
TZ_\zeta\otimes\det f^*TB_\zeta$ be induced from $s_\zeta$ by the sequence above. 
Consider the characteristic class $\Td_g$ as defined in \cite[section 3.3]{KR1}. We
define the map
$\widetilde\Td_g(\mtr{TM},\mtr{TB}):M_g\to\BC$ as
\begin{equation}\label{tdtilde}
\widetilde\Td_g(\mtr{TM},\mtr{TB}):=\Td_g(TM)\sum_{\zeta\neq 1}
\frac1{\zeta-1}\log\frac{\|s_\zeta\|^2}{\|s'_\zeta\|^2}
\end{equation}
(see \cite{Ma2} for the general definition which corresponds to the definition in
\cite[Th. 3.5]{KR1} up to a sign). For our proof, we need
to assume that $f$ extends over Spec $\bf Z$ to a flat projective map of arithmetic
varieties
$f:X\to Y$. We assume furthermore that $g$ corresponds to the action of a diagonalisable
group scheme $\mu_N$ as in \cite[section 4]{KR1} and that $E$ extends to an equivariant
locally free sheaf over $X$.

Consider the structure morphisms $f_X:X\to{\rm Spec}\, {\bf Z}$, $f_Y:Y\to{\rm Spec}\,
{\bf Z}$ and endow the direct image sheafs $R^\cdot f_{X*}E$, $R^\cdot f_*E$ and
$R^\cdot f_{Y*}(R^\cdot f_*E)$ over $\BC$ with the $L^2$-metrics induced by Hodge
theory. The direct image sheaf equipped with these metrics shall be denoted by $R^\cdot
f_{X*}\mtr E$ etc.
The Leray spectral sequence provides an isomorphism
$R^\cdot f_{X*}E\to R^\cdot f_{Y*}(R^\cdot f_*E)$. Consider over $\BC$ the equivariant
Knudsen-Mumford determinant lines $\lambda_{M,g}:=\det_g (H^\cdot(M,\mtr E))^{-1}$ and
$\lambda_{B,g}:=\det_g (\sum_{j,k}(-1)^{j+k}H^j(B,H^k(Z,\mtr E_{|Z})))^{-1}$ and let
$\sigma:\lambda_{B,g}\to\lambda_{M,g}$ denote the isomorphism between them induced by
the Leray spectral sequence. Thus, $\log \|\sigma\|^2_{L^2,g}=\log\frac {\|
s\|^2_{\lambda_M,g}} {\|\sigma s\|^2_{\lambda_B,g}}$ for non-vanishing $s\in
\lambda_{M,g}$.
\begin{theor}\label{MaFormel}
The equivariant analytic torsions of $M$, $B$ and $Z$ are related by the equation
\begin{eqnarray*}
\lefteqn{
\log\|\sigma\|^2_{L^2,g}-T_{g}(M,\mtr E)
}\\&&
=-T_g(B,R^\cdot f_*\mtr E)
-\int_{B^g}\Td_g({TB})T_g(Z,\mtr E)
+\int_{M^g}\widetilde\Td_g(\mtr{TM},\mtr{TB})\ch_g(E)\,\,.
\end{eqnarray*}
\end{theor}
As all fixed points are isolated in our case, it does not make any
difference in this formula whether one considers classes $\Td_g({TB})$ or
characteristic differential forms like
$\Td_g(\mtr {TB})$ as in
\cite{Ma2}.  In \cite{Ma2}, Ma does not need the existence of models over Spec $\bf Z$.
\beginProof We use the notations of \cite[Th. 6.14]{KR1}. Then the arithmetic Lefschetz
formula applied subsequently to
$f_X$, $f$ (in its generalization in
\cite[end of section 2, (b)]{KR2}) and to $f_Y$ yields

\begin{eqnarray*}\lefteqn{
2\deg\ar \ch_{\mu_N}(R^\cdot f_{X*}\mtr E)-T_{g}(M,\mtr E)
=2\deg \ar\ch_{\mu_N}(f_{X!}\mtr E)
}\\&=&
2\deg f^{\mu_N}_{X*}\left(\ar\Td_{\mu_N}(\mtr{Tf_X})\ar\ch_{\mu_N}(\mtr E)\right)
-\int_{M^g}\Td_g(TM)R_g(TM)\ch_g(E)
\\&=&
2\deg
f^{\mu_N}_{Y*}\Bigg[\ar\Td_{\mu_N}(\mtr{Tf_Y})f^{\mu_N}_*\left(
\ar\Td_{\mu_N}(\mtr{Tf})\ar\ch_{\mu_N}(\mtr
E)\right)
\Bigg]
\\&&
+\int_{M^g}\widetilde\Td_g(\mtr{TM},\mtr{TB})\ch_g(E)
-\int_{M^g}\Td_g(TM)R_g(TM)\ch_g( E)
\\&=&
2\deg f^{\mu_N}_{Y*}\Bigg[\ar\Td_{\mu_N}(\mtr{Tf_Y})\Big(
\ar\ch_{\mu_N}(R^\cdot f_* \mtr E)-T_g(Z,\mtr E)
\\&&
+\int_{Z^g}\Td_g(TZ)R_g(TZ)\ch_g(E)
\Big)\Bigg]
\\&&
+\int_{M^g}\widetilde\Td_g(\mtr{TM},\mtr{TB})\ch_g(E)
-\int_{M^g}\Td_g(TM)R_g(TM)\ch_g( E)
\\&=&
2\deg \ar\ch_{\mu_N}(R^\cdot f_{Y*}(R^\cdot f_*\mtr E))-T_g(B,R^\cdot f_*\mtr E)
\\&&
+\int_{B^g}\Td_g(TB)R_g(TB)\ch_g(R^\cdot f_*E)
\\&&-\int_{B^g}\Td_g(\mtr{TB})T_g(Z,\mtr E)
+\int_{B^g}\Td_g(TB)\int_{Z^g}\Td_g(TZ)R_g(TZ)\ch_g( E)
\\&&
+\int_{M^g}\widetilde\Td_g(\mtr{TM},\mtr{TB})\ch_g(\mtr E)
-\int_{M^g}\Td_g(TM)R_g(TM)\ch_g( E)
\\&=&
2\deg \ar\ch_{\mu_N}(R^\cdot f_{Y*}(R^\cdot f_*\mtr E))
-T_g(B,R^\cdot f_*\mtr E)
\\&&
-\int_{B^g}\Td_g(\mtr{TB})T_g(Z,\mtr E)
+\int_{M^g}\widetilde\Td_g(\mtr{TM},\mtr{TB})\ch_g(\mtr E)\,\,.
\end{eqnarray*}
Furthermore
$$
2\deg\ar \ch_{\mu_N}(R^\cdot f_{X*}\mtr E)
-2\deg \ar\ch_{\mu_N}(R^\cdot f_{Y*}(R^\cdot f_*\mtr E))
=\log\|\sigma\|^2_{L^2,g}
$$
by the Leray spectral sequence, which in particular identifies the torsion subgroups
of the direct images. Thus we get the statement of the theorem.
\endProof
{\bf Remark.} By $K$-theory, one could in fact reduce oneself to the case where the
direct images are locally free and non-zero only in degree 0.  Then one can
consider the equivariant determinant lines over Spec $\BZ$ and identify them directly as
graded direct sums of $\BZ$-modules of rank one (up to a sign). Thus one could avoid the
use of the spectral sequence entirely.

An arithmetic
Lefschetz formula for more general fibrations (combined with a general arithmetic
Grothendieck-Riemann-Roch theorem) would provide this formula not only for non-isolated
fixed points, but for the torsion forms associated to a double fibration
$M\to B\to S$ with a K\"ahler manifold $S$ as treated in \cite{Ma1} in the non-equivariant
case.


\section{The torsion of generalized flag manifolds}

Set
 $T_c:=G_c\cap T_\BC$.
We denote the roots of $K$ by $\Sigma_K$; in general, we shall denote the objects
corresponding to subgroups of
$G$ by writing this subgroup as an index. Then
$M:=G_\bC/P_\BC$ is canonically isomorphic to
$G_c/K$ and $\Sigma^+=\Psi\cup \Sigma_K^+$. Recall that we identify ${\frak m}^{1,0}$
with the tangent space $T_e M^{1,0}$ and that the negative of the Killing form induces
an Hermitian metric on $TM$. Fix
$X_0\in{\frak t_c}$ in the closure of the positive Weyl chamber such that its
stabilizer with respect to the adjoint action of
$G_c$ equals
$K$. Then the metric
$(\cdot,\cdot)_{X_0}$ on $\frak m$ which equals for $X\in{\frak m}_\a$, $Y\in{\frak
m}_\beta$
$$
(X,Y)_{X_0}=\left\{
\begin{array}{cc}
0&\alpha\neq\beta\cr
{(X,Y)_B}{ \alpha(X_0)}&\alpha=\beta\cr
\end{array}\right.
$$
induces a K\"ahler metric on $M$ \cite[Ch. 8.D]{Besse}.

Set ${\frak t}_{\rm reg}:=\{X\in{\frak t_c}|\a(X)\notin {\bf Z} \
\forall\a\in\Sigma\}$. For $X\in\frak t_c$ let $e^X\in T_c$ denote the associated group
element. We denote the
$G_c$-representation with highest weight $\lambda$ by $V^G_{\rho+\lambda}$ and its
character is denoted by $\chi_{\rho+\lambda}$. In general, for a weight $\lambda$
and $X\in{\frak t}_{\rm reg}$ we define $\chi_{\rho+\lambda}$ by the Weyl character
formula
$$
\chi_{\rho+\lambda}(e^X):=\frac{\sum_{w\in W_G}(-1)^{l(w)} e^{2\pi i
(\rho+\lambda)(wX)}}{\prod_{\a\in\Sigma^+}2i\sin\pi\a(X)}\,\,.
$$
An irreducible $K$-representation $V^K_{\rho_K+\lambda}$
induces a
$G_c$-invariant holomorphic vector bundle $E^K_{\rho_K+\lambda}$ on $M$. As
$V^K_{\rho_K+\lambda}$ carries a $K$-invariant Hermitian metric which is
unique up to a factor, we get corresponding $G_c$-invariant metrics on
$E^K_{\rho_K+\lambda}$.
Set $\Psi^+:=\{\a\in\Psi|\langle\a^\vee,\rho+\lambda\rangle\geq0\}$ and
$\Psi^-:=\{\a\in\Psi|\langle\a^\vee,\rho+\lambda\rangle<0\}$ with
$\a^\vee=2\a/\|\a\|^2$.

Define for $\phi\in{\bf R}$ and ${\rm Re}\,s>1$
\begin{equation}
\zeta_L(s,\phi)=\sum_{k=1}^\infty \frac{e^{ik\phi}}{k^s}\,\,.
\end{equation}
The function $\zeta_L$ has a meromorphic continuation to the complex plane in $s$ which
is holomorphic for $s\neq 1$. Set $\zeta'_L(s,\phi):=\partial/\partial s
(\zeta_L(s,\phi))$. Let
$P:{\bf Z}\to{\bf C}$ be a function of the form
\begin{equation} P(k)=\sum_{j=0}^m c_jk^{n_j} e^{ik\phi_j}\label{sternnn}
\end{equation} with $m\in{\bf N}_0$, $n_j\in{\bf N}_0$, $c_j\in{\bf C}$,
$\phi_j\in{\bf R}$ for all $j$. We define $P^\odd(k):=(P(k)-P(-k))/2$. Also we define
as in
\cite{K2}
\begin{eqnarray}
\bz P&:=&\sum_{j=0}^m c_j\zeta_L(-n_j,\phi_j),\\
\bzs P&:=&\sum_{j=0}^m c_j\zeta_L'(-n_j,\phi_j),\\
\bzo P&:=&\sum_{j=0}^m c_j\zeta_L(-n_j,\phi_j)\sum_{\ell=1}^{n_j}\frac1
\ell.\\
\mbox{Res }P(p)&:=&\sum_{j=0\atop\phi_j\equiv 0{\rm\ mod }2\pi}^m c_j
\frac{p^{n_j+1}}{2(n_j+1)}\\
\mbox{and}\qquad P^*(p)&:=&-\sum_{j=0\atop\phi_j\equiv 0{\rm\ mod
}2\pi}^m c_j \frac{p^{n_j+1}}{4(n_j+1)} \sum_{\ell=1}^{n_j}\frac 1 \ell
\end{eqnarray}
for $p\in{\bf R}$.
In particular, for $X\in{\frak t}_{\rm reg}$ we get (compare \cite[Th. 10, top of p.
108]{K2})
\begin{equation}\label{zetachi}
\bz\chi_{\rho+\lambda-k\a}(e^X)=\frac{\sum_{w\in W_G}(-1)^{l(w)} e^{2\pi i
(\rho+\lambda)(wX)}(e^{2\pi i\a(wX)}-1)^{-1}}
{\prod_{\beta\in\Sigma^+}2i\sin\pi\beta(X)}\,\,.
\end{equation}
Assume now that $M=G_c/K$ is an Hermitian symmetric space. Then the isotropy
representation is irreducible, hence all $G_c$-invariant metrics coincide with
$(\cdot,\cdot)_B$ up to a factor. Fix one metric associated to
$X_0$ with stabilizer $K$. In \cite[section 11]{K2}, the analytic torsion  on $G_c/K$
has been calculated for vector bundles
$\mtr E^K_{\rho_K+\lambda}$ with $\lambda\in I_+$. We shall now extend this
result to arbitrary $\lambda$.
\begin{prop}
Assume that $G_c/K$ is Hermitian symmetric. Then the zeta function defining the
analytic torsion of $\mtr E^K_{\rho_K+\lambda}$ is given by
\begin{eqnarray*}
Z_t(s)&=&-\sum_{\a\in\Psi^+}\a^\vee(X_0)^s
\sum_{k>\langle\a^\vee,\rho+\lambda\rangle} \frac{\chi_{\rho+\lambda-k
\a}(t)}{k^s(k-\langle\a^\vee,\rho+\lambda\rangle)^s}\\
&&+\sum_{\a\in\Psi^-}\a^\vee(X_0)^s
\sum_{k>-\langle\a^\vee,\rho+\lambda\rangle} \frac{\chi_{\rho+\lambda+k
\a}(t)}{k^s(k+\langle\a^\vee,\rho+\lambda\rangle)^s}\,\,.
\end{eqnarray*}
\end{prop}
\beginProof
According to \cite[(114)]{K2}, the irreducible representations occurring as
eigenspaces of the Laplacian acting on ${\frak A}^{0,q}(M,E)$ are given by the
infinitesimal characters
$\rho+\lambda+k\a$ with
$k>0$.
Consider now the case $\langle\a^\vee,\rho+\lambda\rangle<0$ and assume
$0<k<\langle-\a^\vee,\rho+\lambda\rangle/2$. Let $S_\a$ denote the reflection of the
weights by the hyperplane orthogonal to $\a$. Then the multiplicity of
$V_{\rho+\lambda+k\a}$ cancels with that of
$V_{S_\a(\rho+\lambda+k\a)}=V_{\rho+\lambda+(-k-\langle\a^\vee,\rho+\lambda\rangle)\a}$. In case
$\langle\a^\vee,\rho+\lambda\rangle$ is even, the representation
$V_{\rho+\lambda-\langle\a^\vee,\rho+\lambda\rangle\a/2}$ vanishes. Also, the representation
$V_{\rho+\lambda-\langle\a^\vee,\rho+\lambda\rangle\a}$ corresponds to the eigenvalue 0, thus it
does not contribute to the zeta function. Thus, $Z_t(s)$ is given by
$$
Z_t(s)=\sum_{\a\in\Psi}\a^\vee(X_0)^s
\sum_{k>\max\{0,\langle-\a^\vee,\rho+\lambda\rangle\}} \frac{\chi_{\rho+\lambda+k
\a}}{k^s(k+\langle\a^\vee,\rho+\lambda\rangle)^s}
$$
which equals the above formula after another application of $S_\a$ for $\a\in\Psi^+$.
\endProof

\begin{theor}\label{symm}
Assume that $G_c/K$ is Hermitian symmetric. Then the equivariant analytic torsion of
$\mtr E^K_{\rho_K+\lambda}$ is given by
\begin{eqnarray*}
\lefteqn{
T_t(G_c/K,\mtr E^K_{\rho_K+\lambda})=-2\sum_{\a\in\Psi} \bzs
\chi^\odd_{\rho+\lambda-k\a}
-2\sum_{\a\in\Psi} 
\chi^*_{\rho+\lambda-k\a}(\langle\a^\vee,\rho+\lambda\rangle)}\\
&&-\sum_{\a\in\Psi} \bz
\chi_{\rho+\lambda-k\a}\cdot\log\a^\vee(X_0)
-
\chi_{\rho+\lambda}\sum_{\a\in\Psi^+}\log\a^\vee(X_0)\\
&&-\sum_{\a\in\Psi^+}\sum_{k=1}^{\langle\a^\vee,\rho+\lambda\rangle}
\chi_{\rho+\lambda-k\a}\cdot\log k
+\sum_{\a\in\Psi^-}\sum_{k=1}^{\langle-\a^\vee,\rho+\lambda\rangle}
\chi_{\rho+\lambda+k\a}\cdot\log k\,\,.
\end{eqnarray*}
\end{theor}
\beginProof
According to \cite[Lemma 8]{K2} and \cite[(61),(63)]{K2}, we find
\begin{eqnarray*}
Z_t'(0)&=&\sum_{\a\in\Psi^+}\big[-2\bzs\chi^\odd_{\rho+\lambda-k\a}
-2\chi^*_{\rho+\lambda-k\a}(\langle\a^\vee,\rho+\lambda\rangle)\\
&&-\sum_{k=1}^{\langle\a^\vee,\rho+\lambda\rangle}
\chi_{\rho+\lambda-k\a}\cdot\log k\big]\\
&&
-\sum_{\a\in\Psi^+}\big[
\bz\chi_{\rho+\lambda-k\a}-
\sum_{k=1}^{\langle\a^\vee,\rho+\lambda\rangle} \chi_{\rho+\lambda-k\a}
\big]\log\a^\vee(X_0)\\
&&+\sum_{\a\in\Psi^-}\big[2\bzs\chi^\odd_{\rho+\lambda+k\a}
+2\chi^*_{\rho+\lambda+k\a}(\langle-\a^\vee,\rho+\lambda\rangle)\\
&&+\sum_{k=1}^{\langle-\a^\vee,\rho+\lambda\rangle}
\chi_{\rho+\lambda+k\a}\cdot\log k\big]
-\sum_{\a\in\Psi^-}
\bz\chi_{\rho+\lambda-k\a}
\log\a^\vee(X_0)
\end{eqnarray*}
As $$-\sum_{k=1}^{\langle\a^\vee,\rho+\lambda\rangle} \chi_{\rho+\lambda-k\a}
=-\chi_{\rho+\lambda-\langle\a^\vee,\rho+\lambda\rangle\a}=\chi_{\rho+\lambda}$$ for
$\a\in\Psi^+$ and
$$
\bzs\chi^\odd_{\rho+\lambda+k\a}=-\bzs\chi^\odd_{\rho+\lambda-k\a}\,\,,
$$
$$
\chi^*_{\rho+\lambda+k\a}(\langle-\a^\vee,\rho+\lambda\rangle)=
-\chi^*_{\rho+\lambda-k\a}(\langle\a^\vee,\rho+\lambda\rangle)
$$
we get the statement of the theorem.
\endProof
According to \cite[(103),(116)]{K2}, for any generalized flag manifold the
equation $$\int_{M^t}
\Td_t(TM)R_t(TM)\ch_t(E^K_{\rho_K+\lambda})
=-\sum_{\a\in\Psi}\big[2\bzs \chi^\odd_{\rho+\lambda-k\a}
+\bzo \chi^\odd_{\rho+\lambda-k\a}\big]
$$
holds. For fixed $X\in{\frak t}_{\rm reg}$, according to \cite[eq. (66)]{K2} the
characters $\chi_{\rho+\lambda-k\a}(e^{iX})$ can be written as a linear combination of
exponential functions in
$k$. Polynomial functions in $k$ do not occur in this case. Thus, by definition
$\chi^*_{\rho+\lambda-k\a}$ and
$\bzo
\chi_{\rho+\lambda-k\a}$ vanish.  Hence, we get the following simpler formula for
generic
$X$:
\begin{cor}\label{symmisol}
Assume that $M=G_c/K$ is Hermitian symmetric and choose $X\in{\frak t}_{\rm reg}$. Then
\begin{eqnarray*}
\lefteqn{
T_t((M,g_{X_0}),\mtr E^K_{\rho_K+\lambda})-
\int_{M^t} \Td_t(TM)R_t(TM)\ch_t(E^K_{\rho_K+\lambda})=
}\\&&
-\sum_{\a\in\Psi} \bz
\chi_{\rho+\lambda-k\a}\cdot\log\a^\vee(X_0)
-
\chi_{\rho+\lambda}\sum_{\a\in\Psi^+}\log\a^\vee(X_0)
\\&&
-\sum_{\a\in\Psi^+}\sum_{k=1}^{\langle\a^\vee,\rho+\lambda\rangle}
\chi_{\rho+\lambda-k\a}\cdot\log k
+\sum_{\a\in\Psi^-}\sum_{k=1}^{\langle-\a^\vee,\rho+\lambda\rangle}
\chi_{\rho+\lambda+k\a}\cdot\log k\,\,.
\end{eqnarray*}
\end{cor}

We say that a semisimple Lie group {\bf has tiny weights} if none of its simple
components is of type $G_2$, $F_4$ or $E_8$.
\begin{theor}\label{main}
Let $G_\bC$ be a semisimple Lie group having tiny weights and let $M=G_c/K$ be an
associated (generalized) flag manifold. Let
$\mtr E^K_{\rho_K+\lambda}$ be a
$G_c$-invariant holomorphic Hermitian vector bundle on
$M$ and fix a K\"ahler metric $g_{X_0}$ associated to $X_0\in{\frak t_c}$ on $M$.
Choose
$X\in{\frak t}_{\rm reg}$ and set $t:=e^X$. Then the equivariant analytic torsion
associated to the action of $t$ is given by
\begin{eqnarray*}
\lefteqn{T_t((M,g_{X_0}),\mtr E^K_{\rho_K+\lambda})-
\int_{M^t} \Td_t(TM)R_t(TM)\ch_t(E^K_{\rho_K+\lambda})=
}\\
&&-\sum_{\a\in\Psi} \bz
\chi_{\rho+\lambda-k\a}\cdot\log\a^\vee(X_0)
+
 C \chi_{\rho+\lambda}
\\&&
-\sum_{\a\in\Psi^+}\sum_{k=1}^{\langle\a^\vee,\rho+\lambda\rangle}
\chi_{\rho+\lambda-k\a}\cdot\log k
+\sum_{\a\in\Psi^-}\sum_{k=1}^{\langle-\a^\vee,\rho+\lambda\rangle}
\chi_{\rho+\lambda+k\a}\cdot\log k
\end{eqnarray*}
where the constant $C\in{\bf R}$ does not depend on $t$.
\end{theor}

Our proof proceeds by constructing first a tower of fibrations of flag
manifolds with total space a full flag manifold $G_\bC/B_\BC$, such that all fibers are
Hermitian symmetric spaces. As the torsion is known for the fibers, we can deduce the
value of the torsion for the total space of every fibration from the value on the
base. Thus we finally get the value for the full flag manifold. Then we use the
fibration
$G_\BC/B_\BC\to G_\BC/P_\BC$ with both total space and fiber being full flag manifolds
to compute the torsion for any flag manifold. We need first a Lemma and a
Proposition:

\begin{lemma}\label{typeS}
Let $G_\bC$ be a reductive group having tiny weights. Then there is an ordering
$(\a_1,\dots,\a_m)$ of the base $\Pi$ such that with $\Theta_j:=\{\a_1,\dots,\a_j\}$,
$0\leq j\leq m$, all quotients of subsequent parabolic subgroups
$P_{\Theta_{j+1}}/P_{\Theta_j}$, $0\leq j<m$, are Hermitian symmetric spaces.
\end{lemma}

\beginProof
The only simple groups whose quotients by parabolic subgroups are never Hermitian
symmetric spaces are the groups $G_2$, $F_4$ and $E_8$ \cite[Ch. X \S 6.3]{Hel}. If
$G_\bC$ has tiny weights, then the Dynkin diagram of $\frak g_2$, $\frak f_4$ or
$\frak e_8$ can never occur as a subdiagram of the Dynkin diagram of ${\frak g}_\bC$.
Thus, there is always a simple root
$\a$ such that the quotient $G_\BC/P_{\Pi\setminus\{\a\}}$ is Hermitian symmetric
(namely, those
$\a$ such that their coefficient in the highest root is 1, compare \cite[p.
476ff]{Hel}), and the Levi component of
$P_{\Pi\setminus\{\a\}}$ has again tiny weights. We get the statement by
induction on the number of roots.
\endProof

\begin{prop}\label{special}
Theorem \ref{main} holds for the spaces $G_\bC/{P_{\Theta_j}}$ with $\Theta_j$ as in
Lemma
\ref{typeS}. In particular, it holds for $G_\bC/B_\BC$.
\end{prop}

\beginProof
The proof proceeds by induction on the number of elements of $\Pi\setminus\Theta_j$.
For
$\Theta_j=\Pi$, the analytic torsion vanishes and the statement is trivial. Now assume
that the statement is true for all $\Theta_k\supset\Theta_j$. To shorten the
notation, we shall denote $K_{\Theta_j}$ and $K_{\Theta_{j+1}}$ by $\KO$ and $\K1$,
respectively. Thus, the associated holomorphic fibration
$$
\pi:G_c/\KO\to G_c/\K1
$$
has as fiber the Hermitian symmetric space $\K1/\KO$. The K\"ahler metric on
$G_c/\KO$ induces a K\"ahler metric on $\K1/\KO$. We fix a K\"ahler
metric on $G_c/\K1$ associated to $X_1\in{\frak t_c}$. Notice that $X_1$ and $X_0$
have to be different as their stabilizers are different.

By theorem \ref{ch18} applied to $K/H$ and the base change formulae, we find that a
certain shift of the direct image
$R^\bullet
\pi_* E^{\KO}_{\rho_{\KO}+\lambda}$ is given by the vector bundle
$E^{\K1}_{\rho_{\K1}+\lambda}$ associated to the irreducible
$\K1$-representation with highest weight $\lambda$.
Also, the cohomologies
$H^\bullet(G_c/\KO, E^{\KO}_{\rho_{\KO}+\lambda})$ and
$H^\bullet(G_c/\K1, E^{\K1}_{\rho_{\K1}+\lambda})$ can both be identified with the
irreducible
$G$-representation $E^G_{\rho+\lambda}$, and all
cohomology groups except one vanish. As a $G$-invariant metric on an irreducible
representation is unique up to a constant, the isomorphism $\sigma$ in section
\ref{masection} simply corresponds to a multiplication of the metric with a constant,
and the equivariant metrics
considered in theorem \ref{MaFormel} differ
just by a constant
$C'\in{\bf R}$ times
$\chi_{\rho+\lambda}(t)$.

Let $\widetilde \Td_t(TG_c/\KO,TG_c/\K1,g_{X_0},g_{X_1})$ denote the equivariant
Bott-Chern secondary class associated to the short exact sequence
$$
0\to T\K1/\KO\to TG_c/\KO \to \pi^*TG_c/\K1 \to0
$$
of vector bundles on $G_c/\KO$, equipped with the Hermitian metrics induced by $X_0$,
$X_0$ and $X_1$, respectively.
By theorem \ref{MaFormel} we get the formula
\begin{eqnarray*}
T_t(G_c/\KO,\mtr E^{\KO}_{\rho_{\KO}+\lambda})&=&
\log\|\sigma\|^2_{L^2,t}+T_t(G_c/\K1,\mtr E^{\K1}_{\rho_{\K1}+\lambda})
\\&&
+\int_{(G_c/\K1)^t} \Td_t({TG_c/\K1})
T_t(\K1/\KO,\mtr E^{\KO}_{\rho_{\KO}+\lambda})
\\&&
-\int_{(G_c/\KO)^t} \widetilde \Td_t(TG_c/\KO,TG_c/\K1,g_{X_0},g_{X_1}) \ch_t
( E^{\KO}_{\rho_{\KO}+\lambda})
\end{eqnarray*}
The fixed point set
$(G_c/\K1)^t$ is given by
$W_G/W_\K1$. For simplicity, we assume $\Psi^-=\emptyset$; the proof remains the same
in the general case.

By applying the induction hypothesis and
\cite[Theorem 18]{K2} for the torsion of the fiber, we get
\begin{eqnarray*}
\lefteqn{
T_t(G_c/\KO,\mtr E^{\KO}_{\rho_{\KO}+\lambda})=C'\chi_{\rho+\lambda}(t)+
\int_{(G_c/\K1)^t} \Td_t(TG_c/\K1)R_t(TG_c/\K1)\ch_t(E^{\K1}_{\rho_{\K1}+\lambda})
}\\&&
-\sum_{\alpha\in\Sigma^+\setminus
\Sigma^+_{\K1}}\sum_{k=1}^{\langle\a^\vee,\rho+\lambda\rangle}
\chi_{\rho+\lambda-k\alpha}(t) \log k
-\sum_{\alpha\in\Sigma^+\setminus \Sigma^+_{\K1}}
\bz \chi_{\rho+\lambda-k\alpha}(t)
\log\a^\vee(X_1)
\\&&
+C''\chi_{\rho+\lambda}(t) 
\\&&
+\sum_{w\in W_G/W_\K1} \Td_t(TG_c/\K1)\Big(
\left(\int_{(\K1/\KO)^t}
\Td_t(T\K1/\KO)R_t(T\K1/\KO)\ch_t(E^{\KO}_{\rho_{\KO}+\lambda})\right)_{|wK}
\\&&
-\sum_{\alpha\in\Sigma^+_{\K1}\setminus
\Sigma^+_{\KO}}\sum_{k=1}^{\langle\a^\vee,\rho_{\K1}+\lambda\rangle}
\chi^{\K1}_{\rho_{\K1}+\lambda-k\alpha}(wt) \log k
-\sum_{\alpha\in\Sigma^+_{\K1}\setminus \Sigma^+_{\KO}} 
\bz\chi^{\K1}_{\rho_{\K1}+\lambda-k\alpha}(wt)
\log\a^\vee(X_0)
\\&&
-\chi^\K1_{\rho_\K1+\lambda}(t)\sum_{\alpha\in\Sigma^+_{\K1}\setminus
\Sigma^+_{\KO}}\log\a^\vee(X_0)
\Big)
\\&&
-\int_{(G_c/\KO)^t} \widetilde \Td_t(TG_c/\KO,TG_c/\K1,g_{X_0},g_{X_1}) \ch_t
(\mtr E^{\KO}_{\rho_{\KO}+\lambda})\,\,.
\end{eqnarray*}
Next we notice that
\begin{eqnarray*}
&&\int_{(G_c/\K1)^t} \Td_t(TG_c/\K1)R_t(TG_c/\K1)\ch_t(E^{\K1}_{\rho_{\K1}+\lambda})
\\&&
+\int_{(G_c/\K1)^t} \Td_t(TG_c/\K1)\int_{(\K1/\KO)^t} \Td_t(T\K1/\KO)R_t(T\K1/\KO)\ch_t(E^{\KO}_{\rho_{\KO}+\lambda})
\\
&=&\int_{(G_c/\K1)^t} \Td_t(TG_c/\K1)R_t(TG_c/\K1)
\int_{(\K1/\KO)^t} \Td_t(T\K1/\KO)\ch_t(E^{\KO}_{\rho_{\KO}+\lambda})
\\&&
+\int_{(G_c/\KO)^t} \Td_t(TG_c/\KO)R_t(T\K1/\KO)\ch_t(E^{\KO}_{\rho_{\KO}+\lambda})
\\
&=&\int_{(G_c/\KO)^t} \Td_t(TG_c/\KO)R_t(TG_c/\K1) \ch_t(E^{\KO}_{\rho_{\KO}+\lambda})
\\&&
+\int_{(G_c/\KO)^t} \Td_t(TG_c/\KO)R_t(T\K1/\KO)\ch_t(E^{\KO}_{\rho_{\KO}+\lambda})
\\
&=&\int_{(G_c/\KO)^t} \Td_t(TG_c/\KO)R_t(TG_c/\KO) \ch_t(E^{\KO}_{\rho_{\KO}+\lambda})\,\,.
\end{eqnarray*}
Also, by the classical Lefschetz fixed point formula (in this case already
shown in \cite{Bott}, see also
\cite[Theorem 11, eq. (84)]{K2} and recall (\ref{TdDef})) we get
$$
\sum_{w\in W_G/W_\K1} \Td_t(TG_c/\K1)\chi^{\K1}_{\rho_{\K1}+\lambda-k\alpha}(wt)
=\chi_{\rho+\lambda-k\alpha}(t)\,\,;
$$
furthermore,
$\langle\a^\vee,\rho_{\K1}+\lambda\rangle=\langle\a^\vee,\rho+\lambda\rangle$ for
$\a\in\Sigma_{\K1}$ by
\cite[p. 106]{K2}. By the definition of $\widetilde\Td$
we get
\begin{eqnarray*}
\lefteqn{\int_{(G_c/\KO)^t} \widetilde \Td_t(TG_c/\KO,TG_c/\K1,g_{X_0},g_{X_1}) \ch_t
(E^{\KO}_{\rho_{\KO}+\lambda})}
\\&&=
\sum_{w\in W_G/W_\KO}
\frac{\chi^\KO_{\rho_\KO+\lambda}(wt)}{\prod_{\beta\in\Sigma^+\setminus
\Sigma^+_\KO}(1-e^{-2\pi i
\beta(wt)})}
\sum_{\a\in\Sigma^+\setminus\Sigma^+_K} \frac{\log\alpha(X_0)/\alpha(X_1)}
{e^{2\pi i
\a(wt)}-1}
\\&&=
\sum_{\alpha\in\Sigma^+\setminus \Sigma^+_{\K1}} 
\bz\chi_{\rho+\lambda-k\alpha}(t)
\log\frac{\alpha(X_0)}{\alpha(X_1)}
\,\,.
\end{eqnarray*}
Thus, we find
\begin{eqnarray*}
\lefteqn{
T_t(G_c/\KO,\mtr E^{\KO}_{\rho_{\KO}+\lambda})=(C'+C'')\chi_{\rho+\lambda}(t)+
\int_{(G_c/\KO)^t} \Td_t(TG_c/\KO)R_t(TG_c/\KO) \ch_t(E^{\KO}_{\rho_{\KO}+\lambda})
}\\&&
-\sum_{\alpha\in\Sigma^+\setminus
\Sigma^+_{\K1}}\sum_{k=1}^{\langle\a^\vee,\rho+\lambda\rangle}
\chi_{\rho+\lambda-k\alpha}(t) \log k
-\sum_{\alpha\in\Sigma^+_{\K1}\setminus
\Sigma^+_{\KO}}\sum_{k=1}^{\langle\a^\vee,\rho+\lambda\rangle}
\chi_{\rho+\lambda-k\alpha}(t) \log k
\\&&
-\sum_{\alpha\in\Sigma^+\setminus \Sigma^+_{\K1}} 
\bz\chi_{\rho+\lambda-k\alpha}(t)
\log\a^\vee(X_1)
-\sum_{\alpha\in\Sigma^+_{\K1}\setminus \Sigma^+_{\KO}} 
\bz\chi_{\rho+\lambda-k\alpha}(t)
\log\a^\vee(X_0)
\\&&
-\sum_{\alpha\in\Sigma^+\setminus \Sigma^+_{\K1}} 
\bz\chi_{\rho+\lambda-k\alpha}(t)
\log\frac{\alpha(X_0)}{\alpha(X_1)}
-\chi_{\rho+\lambda}(t)\sum_{\alpha\in\Sigma^+_{\K1}\setminus
\Sigma^+_{\KO}}\log\a^\vee(X_0)
\\&=&
C\chi_{\rho+\lambda}(t)+
\int_{(G_c/\KO)^t} \Td_t(TG_c/\KO)R_t(TG_c/\KO) \ch_t(E^{\KO}_{\rho_{\KO}+\lambda})
\\&&
-\sum_{\alpha\in\Sigma^+\setminus
\Sigma^+_{\KO}}\sum_{k=1}^{\langle\a^\vee,\rho+\lambda\rangle}
\chi_{\rho+\lambda-k\alpha}(t) \log k
-\sum_{\alpha\in\Sigma^+\setminus \Sigma^+_{\KO}} 
\bz\chi_{\rho+\lambda-k\alpha}(t)
\log\a^\vee(X_0)\,\,.
\end{eqnarray*}
\endProof

{\par{\bf Proof of Theorem \ref{main}: }} Consider the fibration $G_c/T_c\to G_c/K$
with fiber $K/T_c$. We equip $G_c/T_c$ by a K\"ahler metric associated to
$X_0\in{\frak t_c}$ with stabilizer $T_c$. Using again
theorem \ref{MaFormel} we get the formula
\begin{eqnarray*}
T_t(G_c/K,\mtr E^K_{\rho_K+\lambda})&=&
-\log|\sigma|^2_{L^2,t}+T_t(G_c/T_c,\mtr E^T_\lambda)
\\&&
-\int_{(G_c/K)_t} \Td_t({TG_c/K})
T_t(K/T_c,\mtr E^T_\lambda)
\\&&
+\int_{(G_c/T)_t} \widetilde \Td_t(TG_c/T_c,TG_c/K,g_{X_0},g_{X_1}) \ch_t
( E^T_\lambda)\,\,.
\end{eqnarray*}
By applying the formulae for $T_t(G_c/T_c,\mtr E^T_\lambda)$ and $T_t(K/T_c,\mtr
E^T_\lambda)$ shown in Proposition
\ref{special}, we get the value of $T_t(G_c/K,\mtr E^K_{\rho_K+\lambda})$ by the
analogue of the calculation above.
\endProof

\section{Description of the $L^2$-metric}

Consider a vector bundle $E^K_{\rho_K+\lambda}$ on $M$ associated to a
$\Sigma^+_K$-dominant weight $\lambda$. In this section we describe the
$L^2$-metric on
$H^*(M,E^K_{\rho_K+\lambda})$ induced by embedding this space into the $C^\infty$
differential forms with coefficients in
$E^K_{\rho_K+\lambda}$ via Hodge theory. As the cohomology is an irreducible
$G$-representation, the metric is uniquely determined by the norm of one element.

Classically \cite{Bott}, $\Gamma^\infty(G_c/K,E^K_{\rho_K+\lambda} \otimes \Lambda^q
T^{*0,1}G_c/K)$ has the $L^2$-dense subspace
\begin{equation}\label{komplex}
\bigoplus_{\mu\in I_+} \Hom_K(V_{\rho+\mu},V^K_{\rho_K+\lambda} \otimes \Lambda^q \Ad)
\otimes V_{\rho+\mu}
\end{equation}
where the section associated to $f\otimes v\in\Hom_K(V_{\rho+\mu},V^K_{\rho_K+\lambda}
\otimes \Lambda^q \Ad)
\otimes V_{\rho+\mu}$ is given by the $K$-equivariant $C^\infty$-function
\begin{eqnarray*}
s:&G_c\to&V^K_{\rho_K+\lambda} \otimes \Lambda^q \Ad\\
&g\mapsto &f(g^{-1}v)\,\,.
\end{eqnarray*}
The direct sum (\ref{komplex}) is a direct sum of complexes. The summand for $\mu$
is the set of $K$-invariants in the tensor product of the complex
$\Hom(V_{\rho+\mu},\Lambda^q \Ad)$ with the space $V^K_{\rho_K+\lambda}
\otimes V_{\rho+\mu}$.
This last complex 
is canonically isomorphic to the standard complex calculating the Lie algebra
cohomology for the nilpotent radical \cite[p. 89, eq. (***)]{Schmidt},
\cite[equation (15.3)]{Bott}. Set $w$ such that $w^{-1}(\rho+\lambda)=:\rho+\lambda_0\in
I_+$. Then the cohomology has either highest weight $\lambda_0$  or
it vanishes in case $\lambda_0\notin I_+$. Assume that $\lambda_0\in I_+$. According to
the proof of Kostant's theorem
\cite[Corollary 3.2.11]{Vogan} the first tensor factor of the summand associated to
$\mu=\lambda_0$ in (\ref{komplex}) consists only of a one-dimensional space in degree
$l(w)$. Thus the cohomology is embedded uniquely in the complex (\ref{komplex}) by the
one-dimensional space
$$
\Hom_K(V_{\rho+\lambda_0},V^K_{\rho_K+\lambda} \otimes \Lambda^{l(w)} \Ad)\,\,.
$$
Now assume that $V^K_{\rho_K+\lambda}$ is one-dimensional and equipped with an
Hermitian metric. Fix one
$v'\in V^K_{\rho_K+\lambda}\otimes \Lambda^{l(w)} \Ad$ of weight $w\lambda_0$ and
norm 1. Choose
$v\in V_{\rho+\lambda_0}$ of weight
$w\lambda_0$ and let $v^\vee\in V^\vee_{\rho+\lambda_0}$ be of weight
$-w\lambda_0$ with
$v^\vee(v)=1$. Then
$v^\vee\otimes v'\in\Hom_K(V_{\rho+\lambda_0},V^K_{\rho_K+\lambda} \otimes
\Lambda^{l(w)}
\Ad)$, hence there is an element of $H^{l(w)}(G_c/K,E^K_{\rho_K+\lambda})$ of
weight
$w\lambda_0$ given by the section
$$
s_0:g\mapsto v^\vee(g^{-1}v)\otimes v'.
$$
which is independent of the choice of $v$.
\begin{lemma}\label{l2norm}
Let $M$ be equipped with the metric associated to $X_0\in{\frak t}_{\rm reg}$. Then
the $L^2$-norm of $s_0$ is given by
$$
|s_0|^2_{L^2}=\prod_{\a\in\Psi}
\frac{\a^\vee(X_0)} {\langle\a^\vee,\rho+\lambda\rangle}\,\,.
$$
\end{lemma}
\beginProof
According to \cite[Lemma 7.23]{BGV}, the computation reduces to that for the
measure $d\mtr g$ on $M$ induced by the Haar measure $dg$ via the formula
$$
|s_0|^2_{L^2}=\prod_{\a\in\Psi} \frac{\a^\vee(X_0)} {\langle\a^\vee,\rho\rangle} 
\int_{G_c/K}|s_0({\mtr g})|^2 d\mtr g
$$
(by \cite[(71)]{K2}, our normalized Haar measure differs from the measure considered in
\cite[Lemma 7.23]{BGV} by the factor in the denominator). By the
$K$-invariance of
$|s_0|^2$, this equals
\begin{eqnarray*}
|s_0|^2_{L^2}&=&\prod_{\a\in\Psi} \frac{\a^\vee(X_0)} {\langle\a^\vee,\rho\rangle}
\cdot
\int_G |v^\vee(g^{-1}v)|^2dg\\
&=&\prod_{\a\in\Psi} \frac{\a^\vee(X_0)} {\langle\a^\vee,\rho\rangle}\cdot({\rm dim}
V_{\rho+\lambda})^{-1}\,\,.
\end{eqnarray*}
The last equality follows by the orthogonality relations (\cite[Th. 4.5]{BtD}). As
$\dim V^K_{\rho_k+\lambda}=1$, the Weyl dimension formula shows
$$
\dim
V_{\rho+\lambda}=\prod_{\a\in\Sigma^+}\frac{\langle\a^\vee,\rho+\lambda\rangle}
{\langle\a^\vee,\rho\rangle}=\prod_{\a\in\Psi}\frac{\langle\a^\vee,\rho+\lambda\rangle}
{\langle\a^\vee,\rho\rangle}\,\,.
$$
Thus the Lemma follows.
\endProof

\section{The main result}

\begin{theor}\label{lines} 
Let $G$ be a Chevalley group having tiny weights and let $P$ be a standard parabolic
subgroup. Consider a $P$-module $A$, free of rank one over $\bf Z$, equipped with an
Hermitian metric on $A_\BC$. Let $\lambda$ be the weight of $A$ and choose $w\in
W_G$ such that
$w^{-1}(\rho+\lambda)=\rho+\lambda_0\in I_+$. We denote the induced line bundle
on the $n+1$-dimensional variety
$X=G/P$ by
$\cal L_\lambda$. 
There is a constant $C\in{\bf R}$ such that the following
identity of formal linear combinations of weights holds
\begin{eqnarray*}
\lefteqn{
-(-1)^{l(w)}\sum_{\mu\in X^*(T)} \mu \log \covol
\mtr{H^{l(w)}(X,{\cal L}_\lambda)}_\mu
}\\&&+\sum_{q=0}^n (-1)^q \sum_{\mu\in  X^*(T)}
\mu \log \# H^q(X,{\cal L}_\lambda)_{\mu,{\rm tor}}
\\&
=&-\frac1{2}\sum_{\a\in\Psi^+}\sum_{k=1}^{\langle\a^\vee,\rho+\lambda\rangle}
\chi_{\rho+\lambda-k\a}\cdot\log k
+\frac1{2}\sum_{\a\in\Psi^-}\sum_{k=1}^{\langle-\a^\vee,\rho+\lambda\rangle}
\chi_{\rho+\lambda+k\a}\cdot\log k\\
&&+\Big(
C-\log\covol \mtr A
\Big)\chi_{\rho+\lambda}\,\,.
\end{eqnarray*}
Furthermore,
$$
C=-\frac1{2}\sum_{\a\in\Psi^+}\log\a^\vee(X_0)
$$
if either $G_c/K$ is Hermitian symmetric or if $\lambda\in I_+$.
\end{theor}

\beginProof
Since $T({\bf C})_{\rm Reg}:=\{t\in T(\BC)_{\rm tor}|\a(t)\neq\beta(t)\ \forall
\a,\beta\in\Sigma, \a\neq\beta\}$ is Zariski-dense in $T(\BC)$, it is enough to prove
the identities of complex numbers obtained by substituting $\mu$ with $\mu(t)\in S^1$
for all
$t\in T(\BC)_{\rm Reg}$. We apply \cite[Th. 6.14]{KR1} for imbeddings
$\mu_N\hookrightarrow T$ satisfying property (Reg) (see definition \ref{propReg})
and for
$t\in\mu_N(\BC)$ being a generator. We adopt the notations from \cite[section 4]{KR1}.

Classically, the free part of the cohomology $H^q(X,{\cal
L}_\lambda)$ vanishes in all degrees except at most for $q=l(w)$ \cite{Bott}.
According to
\cite[Th. 6.14]{KR1},
\begin{eqnarray*}
\lefteqn{
-(-1)^{l(w)}\sum_{\mu\in X^*(\mu_N)} \mu(t) \log\covol \mtr{H^{l(w)}(X,{\cal
L}_\lambda)}_\mu }\\&&
+\sum_{q=0}^n (-1)^q \sum_{\mu \in X^*(\mu_N)} \mu(t) \log \# H^q(X,{\cal
L}_\lambda)_{\mu,{\rm tor}
}\\&=&
\frac1{2}T_t(X_\BC,\mtr L_{\lambda,{\bf C}})-\frac1{2}\int_{X({\bf C})_t}
\Td_t(\mtr {TX}_{\bf C})R_t(\mtr {TX}_{\bf C})\ch_t(\mtr L_\lambda)\\
&&+\widehat{\rm deg} f_*(\widehat\Td_{\mu_N}(\mtr{Tf})\widehat\ch_{\mu_N}(\mtr
L_\lambda))
\,\,.
\end{eqnarray*}
Let $f:X\to{\rm Spec}\,\BZ$ be the structure morphism. As $X^{\mu_N}=X^T$, $T$ acts on
$Tf_{|X^{\mu_N}}$. Let $Tf_\a$ denote the direct summand of $Tf_{|X^{\mu_N}}$ associated
to a weight
$\a$. By Corollary \ref{ch13}, the components of the fixed point scheme are indexed by 
$W_G/W_P$. When restricted to the component
$X^{\mu_N}_w$ for
$[w]\in W_G/W_P$, the equivariant characteristic classes are
given by
$$
\widehat\Td_{\mu_N}(\mtr{Tf})_{|X^{\mu_N}_w}
=\prod_{\a\in\Psi}(1-e^{-2\pi i w\a(t)})^{-1}\cdot
\left(
1+\sum_{\a\in\Psi}\frac{\hat c_1(\mtr{Tf}_{w\a})}{1-e^{2\pi i w\a(t)}}
\right)
$$
(compare the expansion of the equivariant Todd-class in \cite[section 3.3]{KR1}) and
$$
\widehat\ch_{\mu_N}(\mtr {\cal L}_\lambda)_{|X^{\mu_N}_w}
=e^{2\pi i w\lambda(t)}(1+\hat c_1(\mtr {\cal L}_\lambda))\,\,.
$$
Hence we get by Lemma \ref{ch17}
\begin{eqnarray*}
\lefteqn{
\widehat{\rm deg} f_*(\widehat\Td_{\mu_N}(\mtr{Tf})\widehat\ch_{\mu_N}(\mtr
{\cal L}_\lambda))
}\\&=&
\sum_{[w]\in W_G/W_P}\frac{e^{2\pi i w\lambda(t)}}{\prod_{\a\in\Psi}(1-e^{-2\pi i
w\a(t)})}\left[
\sum_{\a\in\Psi} \frac{\widehat{\rm deg} f_*\hat c_1(\mtr{Tf}_{w\a})}{1-e^{2\pi i
w\a(t)}}
+\widehat{\rm deg} f_*\hat c_1(\mtr {\cal L}_\lambda)
\right]
\\&=&\sum_{[w]\in W_G/W_P}\frac{e^{2\pi i
w\lambda(t)}}{\prod_{\a\in\Psi}(1-e^{-2\pi i w\a(t)})}\left[
\frac1{2}\sum_{\a\in\Psi} \frac{-\log\a^\vee(X_0)}{1-e^{2\pi i
w\a(t)}}
-\log \covol \mtr A
\right]
\end{eqnarray*}
Furthermore, \cite[Th. 11, eq. (84)]{K2} shows
$$
\sum_{[w]\in W_G/W_P}\frac{-e^{2\pi i
w\lambda(t)}\log \covol \mtr A }{\prod_{\a\in\Psi}(1-e^{-2\pi i w\a(t)})}
=-\chi_{\rho+\lambda}(t)\cdot \log\covol\mtr A
$$
Similarly we shall proceed as
in the proof of \cite[theorem 10]{K2} to obtain a formula for the other term. For a
class $[w]\in W_G/W_P$ we shall denote by $l(w)$ the minimal length of its
representents.
\begin{eqnarray*}
\lefteqn{
\sum_{[w]\in W_G/W_P}\frac{e^{2\pi i
w\lambda(t)}}{\prod_{\a\in\Psi}(1-e^{-2\pi i w\a(t)})}
\sum_{\a\in\Psi} \frac{-\log\a^\vee(X_0)}{1-e^{2\pi i
w\a(t)}}
}\\&=&
\sum_{[w]\in W_G/W_P} \frac{e^{2\pi i w(\rho-\rho_K+\lambda)(t)}(-1)^{l(w)} \sum_{w'\in
W_P} e^{2\pi iww'\rho_K(t)}}{\prod_{\beta\in\Sigma^+}2i \sin\pi\beta(t)}
\sum_{\a\in\Psi}\frac{\log\a^\vee(X_0)}{e^{2\pi i w\a(t)}-1}
\\&=&
\sum_{[w]\in W_G/W_P} \frac1{\prod_{\beta\in\Sigma^+}2i
\sin\pi\beta(t)}
\\&&\cdot
\sum_{\a\in\Psi}\sum_{w'\in W_P} (-1)^{l(w)+l(w')} e^{2\pi
i ww'(\rho+\lambda)(t)} 
 \frac {\log(w'\a)^\vee(X_0)}{e^{2\pi i ww'\a(t)}-1}
\end{eqnarray*}
(in the last equation we used the facts that $\lambda$ and $\Psi$ are
$W_P$-invariant)
$$
=\sum_{\a\in\Psi} \sum_{w\in W_G} \frac{(-1)^{l(w)}e^{2\pi
iw(\rho+\lambda)(t)}(e^{2\pi iw\a(t)}-1)^{-1}} {\prod_{\beta\in\Sigma^+}2i
\sin\pi\beta(t)}
\log \a^\vee(X_0)
$$
(as $X_0$ is $K$-stable)
$$
=
\sum_{\a\in\Psi} \bz \chi_{\rho+\lambda-k\a}(t)\cdot\log\a^\vee(X_0)
$$
(by equation (\ref{zetachi})). Thus
\begin{eqnarray*}
\lefteqn{
\widehat{\rm deg} f_*(\widehat\Td_{\mu_N}(\mtr{Tf})\widehat\ch_{\mu_N}(\mtr
{\cal L}_\lambda))
}\\&=&
\frac1{2} \sum_{\a\in\Psi} \bz
\chi_{\rho+\lambda-k\a}(t)\cdot\log\a^\vee(X_0)-\chi_{\rho+\lambda}(t)\cdot
\log\covol\mtr A\,\,.
\end{eqnarray*}
The proof is finished by combining this result with Th. \ref{main}. For
Hermitian symmetric spaces, combine instead with the more precise corollary
\ref{symmisol} which does not involve the unknown constant $C$. Now consider the
case $\lambda\in I_+$. In this case, by the Kempf
vanishing theorem
\cite[Prop. 4.5]{Jan} the cohomology has no torsion. Furthermore $\Psi=\Psi^+$, thus
the formula in theorem
\ref{lines} simplifies to
\begin{eqnarray*}
\lefteqn{
-\sum_{\mu\in X^*(T)} \mu \log \covol
\mtr{H^0(X,{\cal L}_\lambda)}_\mu
}
\\&
=&-\frac1{2}\sum_{\a\in\Psi}\sum_{k=1}^{\langle\a^\vee,\rho+\lambda\rangle}
\chi_{\rho+\lambda-k\a}\cdot\log k
+\Big(
C-\log\covol \mtr A
\Big)\chi_{\rho+\lambda}\,\,.
\end{eqnarray*}
In particular, the component associated to $\mu=\lambda$ verifies the equality
$$
-\log \covol
\mtr{H^0(X,{\cal L}_\lambda)}_\lambda
=\frac1{2}\sum_{\a\in\Psi}
\log \langle\a^\vee,\rho+\lambda\rangle
+\Big(
C-\log\covol \mtr A
\Big)\,\,,
$$
as the weight $\lambda$ occurs only in the characters $\chi_{\rho+\lambda}$
and $\chi_{\rho+\lambda-\langle\a^\vee,\rho+\lambda\rangle\a}$, with multiplicity
$1$ and $-1$, respectively.
By combining Proposition
\ref{ch110} and Lemma
\ref{l2norm}, we find on the other hand
$$
-\log \covol
\mtr{H^0(X,{\cal L}_\lambda)}_\lambda
=-\frac1{2}\sum_{\a\in\Psi}\log\a^\vee(X_0)
+\frac1{2}\sum_{\a\in\Psi}
\log \langle\a^\vee,\rho+\lambda\rangle
-\log\covol \mtr A
$$
(in Lemma \ref{l2norm}, $\mtr A$ was normalized to have covolume 1).
Thus we get the value of the constant $C$ for dominant $\lambda$.
\endProof

\begin{cor}\label{cor72}
Let $A_{\rho+\lambda}:=H^{l(w)}(X,{\cal L}_\lambda)_{\rm free}$ be a $G$-module induced
by a line bundle
${\cal L}_\lambda$. If $A_{\rho+\lambda}\neq 0$ (i.e. $\lambda_0\in I_+$) equip
$A_{\rho+\lambda}$ with the unique
$G_c$-invariant metric such that the generator of the weight space to the highest
weight $\lambda_0$ has norm 1. Then
\begin{eqnarray*}
\lefteqn{
-(-1)^{l(w)}\sum_{\mu\in X^*(T)} \mu \log \covol \mtr A_{{\rho+\lambda},\mu}
}\\&&
+\sum_{q=0}^n (-1)^q \Big(\sum_{\mu\in X^*(T)} \mu \log \# H^q(X,{\cal L}_\lambda)_{\mu,{\rm
tor}}-(-1)^{l(w)}\chi_{\rho+\lambda}\log \# H^q(X,{\cal L}_\lambda)_{\lambda_0,{\rm
tor}}\Big)
\\&
=&-\frac1{2}\sum_{\a\in\Psi^+}\sum_{k=1}^{\langle\a^\vee,\rho+\lambda\rangle-1}
\chi_{\rho+\lambda-k\a}\cdot\log k
+\frac1{2}\sum_{\a\in\Psi^-}\sum_{k=1}^{\langle-\a^\vee,\rho+\lambda\rangle-1}
\chi_{\rho+\lambda+k\a}\cdot\log k
\end{eqnarray*}
\end{cor}
\beginProof
This follows from theorem \ref{lines} by comparing the components of weight
$\lambda_0$.  The weight $\lambda_0$ occurs only in the characters $\chi_{\rho+\lambda}$
and $\chi_{\rho+\lambda-\langle\a^\vee,\rho+\lambda\rangle\a}$, with multiplicity
$(-1)^{l(w)}$ and $-(-1)^{l(w)}$, respectively.
\endProof
{\bf Remark.} Let $P'\subseteq P$ be another standard parabolic subgroup. Then $A$ is
a $P'$-representation, and both sides of the above formula remain the same. Namely,
$H^\bullet(G/P,{\cal L}_\lambda)=H^\bullet(G/P',{\cal L}_\lambda)$ and
$\chi_{\rho+\lambda-k\a}=0$ for any $\a\in\Sigma^+_K$, $1\leq k\leq
\langle\a^\vee,\rho+\lambda\rangle-1$.

For two $\BZ$-free $G$-modules $A$, $A'$ with $A_\BQ$ and $A'_{\bf
Q}$ isomorphic and irreducible we denote by $[A,A']\in{\bf Q}^+$ the index obtained by
embedding
$A'$ in such a way in
$A\otimes{\bf Q}$ that the weight spaces  of highest weight are identified. Then
$[A,A']=\covol A'/\covol A$. Let $w_0$ be the Weyl group element of maximal
length. When applied to
$w_0\ldotp\lambda$ instead of $\lambda$, the above corollary reads
\begin{eqnarray*}
\lefteqn{
-\frac1{2}\sum_{\a\in\Psi^+}\sum_{k=1}^{\langle\a^\vee,\rho+\lambda\rangle-1}
\chi_{\rho+\lambda-k\a}\cdot\log k
+\frac1{2}\sum_{\a\in\Psi^-}\sum_{k=1}^{\langle-\a^\vee,\rho+\lambda\rangle-1}
\chi_{\rho+\lambda+k\a}\cdot\log k
}\\&
=&-\frac{(-1)^{l(w_0)}}2\sum_{\a\in\Psi^+}\sum_{k=1}^{-\langle-w_0\a^\vee,
w_0(\rho+\lambda)\rangle-1}
\chi_{w_0(\rho+\lambda-k\a)}\cdot\log k
\\&&
+\frac{(-1)^{l(w_0)}}2\sum_{\a\in\Psi^-}\sum_{k=1}^{\langle-w_0\a^\vee,
w_0(\rho+\lambda)\rangle-1}
\chi_{w_0(\rho+\lambda+k\a)}\cdot\log k
\\&
=&-\frac{(-1)^{l(w_0)}}2\sum_{\a\in-w_0\Psi^+}\sum_{k=1}^{-\langle\a^\vee,
w_0(\rho+\lambda)\rangle-1}
\chi_{w_0(\rho+\lambda)+k\a}\cdot\log k
\\&&
+\frac{(-1)^{l(w_0)}}2\sum_{\a\in-w_0\Psi^-}\sum_{k=1}^{\langle\a^\vee,
w_0(\rho+\lambda)\rangle-1}
\chi_{w_0(\rho+\lambda)-k\a}\cdot\log k
\\
&=&
(-1)^{l(w)}\sum_{\mu\in X^*(T)} \mu \log \covol \mtr A_{{w_0(\rho+\lambda)},\mu}
\\&&
-\sum_{q=0}^n (-1)^{q+l(w_0)} \Big(\sum_{\mu\in X^*(T)} \mu \log \#
H^q(X,{\cal L}_{w_0\ldotp\lambda})_{\mu,{\rm tor}}
\\&&
-(-1)^{l(w_0w)}\chi_{w_0(\rho+\lambda)}\log
\# H^q(X,{\cal L}_{w_0\ldotp\lambda})_{{\lambda_0},{\rm tor}}\Big)
\\
&=&
(-1)^{l(w)}\sum_{\mu\in X^*(T)} \mu \log \covol \mtr A_{{w_0(\rho+\lambda)},\mu}
\\&&
-\sum_{q=0}^n (-1)^{q} \Big(\sum_{\mu\in X^*(T)} \mu \log \#
H^{n-q}(X,{\cal L}_{w_0\ldotp\lambda})_{\mu,{\rm tor}}
\\&&
-(-1)^{l(w)}\chi_{\rho+\lambda}\log
\# H^{n-q}(X,{\cal L}_{w_0\ldotp\lambda})_{{\lambda_0},{\rm tor}}\Big)
\,\,.
\end{eqnarray*}
By subtracting the
equations for
$A_{\rho+\lambda}$ and
$A_{w_0(\rho+\lambda)}$  we get an identity which is an
immediate consequence of Corollary
\ref{Serretors} and Proposition \ref{ch19}. Thus corollary \ref{cor72}
is equivalent to the statement obtained by adding  the
equations for
$A_{\rho+\lambda}$ and
$A_{w_0(\rho+\lambda)}$.
We get
\begin{cor}(Jantzen sum formula \cite[p. 311]{Jan})
\begin{eqnarray*}
\lefteqn{
(-1)^{l(w)}\sum_{\mu\in X^*(T)} \mu \log [A_{\rho+\lambda,\mu}:A_{{w_0(\rho+\lambda)},\mu}]
}\\&&
+\sum_{\mu\in X^*(T)} \sum_{q=0}^n \mu\cdot (-1)^q  \Big(\log \# H^q(X,{\cal
L}_\lambda)_{\mu,{\rm tor}}
-\log \#
H^{n-q}(X,{\cal L}_{w_0\ldotp\lambda})_{\mu,{\rm tor}}\Big)
\\&&
-(-1)^{l(w)}\chi_{\rho+\lambda}\sum_{q=0}^n (-1)^q \Big(\log \#
H^q(X,{\cal L}_\lambda)_{\lambda_0,{\rm tor}}-\log
\# H^{n-q}(X,{\cal L}_{w_0\ldotp\lambda})_{{\lambda_0},{\rm tor}}\Big)
\\&
=&-\sum_{\a\in\Psi^+}\sum_{k=1}^{\langle\a^\vee,\rho+\lambda\rangle-1}
\chi_{\rho+\lambda-k\a}\cdot\log k
+\sum_{\a\in\Psi^-}\sum_{k=1}^{\langle-\a^\vee,\rho+\lambda\rangle-1}
\chi_{\rho+\lambda+k\a}\cdot\log k\,\,.
\end{eqnarray*}
\end{cor}
As above, both sides are independent of the parabolic subgroup. In particular, for
$P=B$ we get the usual formulation of the Jantzen sum formula. 

Combining these formulae with \cite[(8.7)]{Jan}
$$
\#
H^q(X,{\cal L}_\lambda)_{\lambda_0,{\rm tor}}=0\qquad\forall\lambda,q
$$
one notices the vanishing of the multiples of $\chi_{\rho+\lambda}$ in the above two
corollaries. We did not use this result before to point out what exactly our result
says without applying further representation theory. For $\lambda\in I_+$, the Kempf
vanishing theorem
\cite[Prop. 4.5]{Jan} states that the torsion of the cohomology vanishes.  Note
 that the other way round the Jantzen sum formula provides the values of the equivariant
Ray-Singer torsion for isolated fixed points up to a multiple of $\chi_{\rho+\lambda}$,
in particular for the types
$G_2$,
$F_4$ and
$E_8$.

\section{The height of generalized flag varieties}

For the definitions of the objects in Arakelov geometry which are used in this
chapter we refer the reader to
\cite{Soule}. According to \cite[14.7]{BH}, the very ample line bundles on $G_c/K$
are induced by the $K$-representations with highest weight $\lambda$ such that
$$
\langle\a^\vee,\lambda\rangle=\left\{
\begin{array}{cc}
=0& \mbox{if }\a\in\Sigma_K\\
>0& \mbox{if }\a\in\Psi\,\,.
\end{array}\right.
$$
Equip a line bundle ${\cal L}_\lambda\to G/P$ to such a $\lambda$ with the
equivariant metric induced by normalizing the length of the generator of the
$P$-module to 1. Thus, as the equivariant metric on ${\cal L}_\lambda$ is unique up to
a factor, it coincides with the metric induced by the ${\cal O}(1)$ bundle via the
projective embedding associated to ${\cal L}_\lambda$. In particular, this metric is
positive (this can be shown directly, too). In this section we compute the {\bf
global height}  of
$G/P$ with respect to
$\mtr {\cal L}_\lambda$ defined as
$$
h(G/P,\mtr {\cal L}_\lambda):=f_*\left[\hat c_1(\mtr
{\cal L}_\lambda)^{n+1}\right]\,\,.
$$
{\bf Remark.} The line bundle ${\cal L}_\lambda$ is also very ample
\cite[II.8.5]{Jan} (we shall not need this fact). See \cite{SZh}, \cite{Abbes} for
relations between the global height of a variety and the height of points on that
variety in this case.

\bigskip
Set $\Psi_j:=\{\a\in\Psi|\langle\a^\vee,\lambda\rangle=j\}$ for
$j\in{\bf N}$. These sets are
$W_K$-invariant as $\Psi$ and $\lambda$ are $W_K$-invariant. Hence for every $j$
there is a virtual $P$-representation with character equal to
$\sum_{\a\in\Psi_j}e^{2\pi i\a}$, defining a virtual vector bundle $E_j$ on
$G/P$. As before, we equip $G_c/K$ with the metric associated to some $X_0$ in a
certain subset of ${\frak t}_c$.

Let Ht denote the additive topological characteristic class associated to the power
series
$$
{\rm Ht}(x):=\sum_{k=0}^\infty \frac{(-x)^k}{2(k+1)(k+1)!}\,\,.
$$
This is the Taylor expansion of the function
$x\mapsto\frac1{2x}\left(\log|x|-\Gamma'(1)-{\rm Ei}(-x)\right)$ at $x=0$ with Ei
being the exponential-integral function \cite[39 (13)]{N}. For $x<0$ ${\rm Ei}(x)$ is
given by
$$
{\rm Ei}(x)=\int_{-\infty}^x\frac{e^t}t\,dt\,\,.
$$
\begin{theor}\label{height}
The height of $G/P$ with respect to $\mtr
{\cal L}_\lambda$ is given by
$$
h(G/P,\mtr {\cal L}_\lambda)=(n+1)!\int_M
{\rm Ht}(\sum_{j\in{\bf N}} j\psi^j E_j) e^{c_1(L_\lambda)}
$$
where $\psi^j$ denotes the $j$-th Adams operator.
\end{theor}

In particular the height can be written as
\begin{equation}\label{harmo}
\sum_{l=0}^n \frac{(-1)^l}{2(l+1)}{n+1\choose l+1}\int_M l!\sum_{j\in{\bf N}}
j^{l+1}{\rm ch}(E_j)^{[l]}c_1(L_\lambda)^{n-l}\,\,.
\end{equation}
We conclude that $h(G/P,\mtr {\cal L}_\lambda)$ has the form
$$
\sum_{\ell=1}^{n+1} \frac{k_\ell}{2\ell}
$$
with $k_\ell\in{\bf Z}\ \forall\ell$. See Corollary \ref{hdim} for a refinement and a
conjecture.

\beginProof
 As in \cite[VIII.2.3]{Soule} we have for $m\to\infty$
\begin{equation}\label{degh}
r(m):=\widehat{\rm deg}(\pi_*\mtr
{\cal L}_{m\lambda},\|\cdot\|^2_{L^2})=
\frac{m^{n+1}}{(n+1)!} h(G/P,\mtr {\cal L}_\lambda)
+O(m^n\log m)\,\,.
\end{equation}
Here we use the fact that $\mtr {\cal L}_\lambda$ is  positive, which implies  a
result by Bismut-Vasserot \cite{BVa} on the asymptotics of the non-equivariant 
holomorphic torsion. Consider the polynomials
$d_{\rho+m\lambda-k\a}:=\dim V_{\rho+m\lambda-k\a}$
$(\a\in\Psi)$ in $m$ and $k$, given by the Weyl dimension
formula
$$
d_{\rho+m\lambda-k\a}=\prod_{\beta\in\Sigma^+}\left(1+{\langle\beta^\vee,m\lambda-k\a\rangle\over
\langle\beta^\vee,\rho\rangle}\right)\,\,.
$$
The term $r(m)$ equals the left hand side of theorem \ref{lines} evaluated at zero.
By the Jantzen sum formula (or theorem
\ref{lines} for
$G$ having tiny weights), for $m\geq1$
\begin{eqnarray*}
r(m)&=&-\frac1{2}\sum_{\a\in\Psi}\sum_{k=1}^{\langle\a^\vee,\rho+m\lambda\rangle}
d_{\rho+m\lambda-k\a}\log
k\\&&
-d_{\rho+m\lambda}\cdot\Big[\log \covol \mtr{H^0(X,{\cal L}_{m\lambda})}_{{m\lambda}}
+\frac1{2}\sum_{\a\in\Psi}\log\langle\a^\vee,\rho+m\lambda\rangle \Big]
\,\,.
\end{eqnarray*}
Here we replaced the constant $C$ in theorem \ref{lines} by the value one obtains by
comparing the components of weight $m\lambda$, similar to the last section in the proof
of theorem
\ref{lines}. As in the case $\lambda\in I_+$ in theorem \ref{lines} (where only groups
$G$ with tiny weights were allowed), we find again by combining Proposition \ref{ch110}
and Lemma
\ref{l2norm}
$$
\log \covol \mtr{H^0(X,{\cal L}_{m\lambda})}_{{m\lambda}}
+\frac1{2}\sum_{\a\in\Psi}\log\langle\a^\vee,\rho+m\lambda\rangle
=\frac1{2}\sum_{\a\in\Psi}\log \a^\vee(X_0)\,\,,
$$
thus this factor is independent of $m$. The expression for $r(m)$ is a sum over terms
which (at a first sight) look like having order $O(m^{1+\#\Sigma^+}\log m)$ for
$m\to\infty$. Thus, it is not obvious that it is in fact of order $O(m^{n+1})$, and
one has to be careful when calculating the term of highest degree.
 We shall need the following
three facts about the polynomials
$d_{\rho+m\lambda-k\a}$:
\begin{enumerate}
\item The degree in $m$ of $d_{\rho+m\lambda-k\a}$ is equal to $n$, since
$$
d_{\rho+m\lambda-k\a}=\prod_{\beta\in\Psi}\frac{\langle\beta^\vee,
\rho+m\lambda-k\a\rangle} {\langle\beta^\vee,
\rho\rangle}\cdot
\prod_{\beta\in\Sigma_K^+}\frac{\langle\beta^\vee,
\rho-k\a\rangle} {\langle\beta^\vee,
\rho\rangle}\,\,.
$$
\item For all $j\in{\bf N}$, the common polynomial degree in $m$ and $k$ of
$\sum_{\a\in\Psi_j} d_{\rho+m\lambda-k\a}$ is less or equal to $n$. This is a
consequence of the Riemann-Roch theorem, which states in this case
$$
\sum_{\a\in\Psi_j} d_{\rho+m\lambda-k\a}
=\int_M \Td(TM)\ch(\psi^{-k}E_j)e^{m c_1(L_\lambda)}
$$
(compare \cite[Th. 13]{K2}).
\item The polynomial $d_{\rho+m\lambda-k\a}$ is skew-symmetric in $k$ around
$k=\langle\a^\vee,\rho+m\lambda\rangle/2$. More precisely, by applying the reflection
$S_\a$ at the hyperplane orthogonal to $\a$ we get
$$
d_{\rho+m\lambda-(k+\langle\a^\vee,\rho+m\lambda\rangle)\a}=-d_{\rho+m\lambda+k\a}\,\,.
$$
\end{enumerate}
Now we see that
\begin{eqnarray*}
-2r(m)+O(m^n)&=&
\sum_j\sum_{\a\in\Psi_j}\sum_{k=1}^{j m}
d_{\rho+m\lambda-k\a}\log k
\\
&&+\sum_{\a\in\Psi}\sum_{k=\langle\a^\vee,m\lambda\rangle+1}^{\langle\a^\vee,\rho+m\lambda\rangle}
d_{\rho+m\lambda-k\a}\log k\\
&=&\sum_j\sum_{k=1}^{j m}\sum_{\a\in\Psi_j}
d_{\rho+m\lambda-k\a}\log k\\
&&-\sum_{\a\in\Psi}\sum_{k=1-\langle\a^\vee,\rho\rangle}^0
d_{\rho+m\lambda+k\a}\log (k+\langle\a^\vee,\rho+m\lambda\rangle)\\
\end{eqnarray*}
using the skew-symmetry in the last equation. Notice that the last double sum is of
order $O(m^n\log m)$. By approximating the first sum with integrals and doing partial
integration we get
\begin{eqnarray*}
-2r(m)&=&\sum_j\int_1^{j m}\sum_{\a\in\Psi_j}
d_{\rho+m\lambda-k\a}\log k \,dk
+O(m^n\log m)\\
&=&\sum_{\a\in\Psi}\int_1^{\langle\a^\vee,m\lambda\rangle}
d_{\rho+m\lambda-k\a}\log k \,dk+O(m^n\log m)\\
&=&\sum_{\a\in\Psi}\int_0^x
d_{\rho+m\lambda-k\a}dk\log x\Big|_{x=1}^{\langle\a^\vee,m\lambda\rangle}\\
&&-\sum_{\a\in\Psi}\int_1^{\langle\a^\vee,m\lambda\rangle}\int_0^x
d_{\rho+m\lambda-k\a}dk\frac{dx}x+O(m^n\log m)\\
&=&\sum_{\a\in\Psi}\int_0^{\langle\a^\vee,m\lambda\rangle}
d_{\rho+m\lambda-k\a}dk\log
{\langle\a^\vee,m\lambda\rangle}\\
&&-\sum_{\a\in\Psi}\int_0^{\langle\a^\vee,m\lambda\rangle}\int_0^x
d_{\rho+m\lambda-k\a}dk\frac{dx}x+O(m^n\log m)\\
\end{eqnarray*}
By the skew-symmetry of $d_{\rho+m\lambda-k\a}$, this equals
\begin{eqnarray}\nonumber
&&\sum_{\a\in\Psi}\int_0^{\langle\a^\vee,\rho\rangle}
d_{\rho+m\lambda-k\a}dk\log
{\langle\a^\vee,m\lambda\rangle}\\ \nonumber
&&-\sum_{\a\in\Psi}\int_0^{\langle\a^\vee,m\lambda\rangle}\int_0^x
d_{\rho+m\lambda-k\a}dk\frac{dx}x+O(m^n\log m)\\
&=&-\sum_{\a\in\Psi}\int_0^{\langle\a^\vee,m\lambda\rangle}\int_0^x
d_{\rho+m\lambda-k\a}dk\frac{dx}x+O(m^n\log m)\label{hopp}
\end{eqnarray}

In particular, we see that $r(m)$ is
of order $O(m^{n+1})$, which can also be derived from equation (\ref{degh}). By
decomposing into the $\Psi_j$-parts again and applying the Riemann-Roch theorem we get
\begin{eqnarray*}
2r(m)&=&\sum_j\int_0^{j m}\int_0^x \int_M \Td(TM){\rm
ch}(\psi^{-k}E_j)e^{mc_1({\cal L}_\lambda)}dk\frac{dx}x+O(m^n\log m)\\
&=&\int_0^m\int_0^x \int_M \Td(TM){\rm
ch}(\sum_j j\psi^{-j k}E_j)e^{mc_1({\cal L}_\lambda)}dk\frac{dx}x+O(m^n\log
m)\\ &=&\int_0^m\int_0^x \int_M {\rm
ch}(\sum_j j\psi^{-j k}E_j)e^{mc_1({\cal L}_\lambda)}dk\frac{dx}x+O(m^n\log m)
\end{eqnarray*}
Hence we can express $r(m)$ as
$$
r(m)=\int_M {\rm
Ht}(\sum_j j\psi^j E_j)e^{c_1({\cal L}_\lambda)}\cdot m^{n+1}+O(m^n\log m)\,\,.
$$
This proves the above theorem.
\endProof
\par{\bf Remark.} If $G_c/K$ is Hermitian symmetric a formula for $T(G_c/K,\mtr L)$ is
known \cite{K2}. In this case we could have equally well worked with the difference
$T_g(G_c/K,\mtr L)-T(G_c/K,\mtr L)$ by arguing similarly as above, but avoiding the
use of
\cite{BVa}. In the case of Hermitian symmetric spaces there is a unique primitive
$\lambda$. Using the classification of irreducible Hermitian symmetric spaces
\cite{Hel}, one verifies that
$\langle\a^\vee,\lambda\rangle$ equals either 1 or 2 for all $\a\in\Psi$.

\medskip
Equation (\ref{hopp}) provides a very effective way to compute the
height. Namely, for $j\in{\bf N}$ consider the sum
\begin{eqnarray}\label{hupp}
f_j(m,k)&:=&\sum_{\alpha\in\Psi_j} d_{\rho+m\lambda-k\alpha}\\
&=&\sum_{\alpha\in\Psi_j}
\prod_{\beta\in\Sigma^+_K} (1-k\frac{\langle\beta^\vee,\alpha\rangle}
{\langle\beta^\vee,\rho\rangle})
\prod_{j'}\prod_{\beta\in\Psi_{j'}} (1+\frac{j'm-k\langle\beta^\vee,\alpha\rangle}
{\langle\beta^\vee,\rho\rangle})\nonumber
\end{eqnarray}
and replace every power $k^l$ by $\frac{(m j)^{l+1}}{2(l+1)^2}$. The height is
obtained by adding the coefficients of $m^{n+1}$ and multiplying with $(n+1)!$. 
For $\a\in\Sigma$
\begin{equation}\label{cox}
\#\{\beta\in\Sigma^+|\langle\a,\beta\rangle\neq 0\}\leq 2 c(G)-3
\end{equation}
where $c(G)$ is the Coxeter number of $G$ \cite[Ch. V \S6.1]{Bour}. Namely, by 
\cite[Ch. VI \S1, Prop. 32]{Bour}, equality holds in (\ref{cox}) for root systems
where all roots have the same length. Furthermore, for $B_l$ and $C_l$ the
cardinality in (\ref{cox}) equals
$2c(G)-5$ and
$c(G)-1$ (depending on the root), for $F_4$ it equals $2c(G)-9=15$ and for $G_2$ it
equals
$c(G)-1=5$. Thus the degree in
$k$ of
$f_j(m,k)$ is less or equal to $2 c(G)-3$, and we find by using the same
reformulation as in formula (\ref{harmo})
\begin{cor}\label{hdim}
There are $k_\ell\in{\bf Z}\ (1\leq \ell\leq 2 c(G)-2)$ such that
$$
h(G/P, \mtr {\cal L}_\lambda)=\sum_{\ell=1}^{2 c(G)-2}
\frac{k_\ell}{2\ell}\,\,.
$$
In other words, the largest power of a prime occurring in the denominator of $2 h(G/P,
\mtr {\cal L}_\lambda)$ is less or equal to $2c(G)-2$.
\end{cor}
Our computations of examples as well as Tamvakis' results in \cite{T1}, \cite{T2},
\cite{T3} strongly suggest the

{\bf Conjecture }{\it There are $k_\ell\in{\bf Z}\ (1\leq \ell\leq c(G)-1)$ such that}
$$
h(G/P, \mtr {\cal L}_\lambda)=\sum_{\ell=1}^{c(G)-1}
\frac{k_\ell}{2\ell}\,\,.
$$

Also there is the following fixed point
expression for $h(G/P, \mtr {\cal L}_{\lambda})$:
\begin{lemma}\label{hfixpt}
Let $Y\in{\frak g}$ act with isolated fixed points and let $\phi$, $\theta_\nu$ denote
the angles of the action on $L_{\lambda|p}$, $TM_{|p}$ respectively
for $p\in M^Y$. Then
$$
h(G/P,\mtr {\cal L}_\lambda)=\sum_{p\in
M^Y}\frac{1}{\prod_\nu\theta_\nu}\sum_{l=1}^{n+1}
\sum_{j\in\BN}\sum_{\theta_\nu\in\Psi_j}
\frac{\phi^{n+1}-\phi^{n+1-l}(\phi-j\theta_\nu)^l}{2l\theta_\nu}\,\,.
$$
\end{lemma}
Note that if $Y$ is an element of the Lie algebra of the maximal torus then there is a
canonical isomorphism
$M^Y=W_G/W_K$. The angle $\phi$ at $[w]\in W_G/W_K$ is given by $(2\pi w\lambda)(Y)$
and the angle $\theta_{w\alpha}$ ($\alpha\in\Psi$) corresponding to $TM_{w\alpha}$ is
given by
$(2\pi w \alpha)(Y)$. Thus, the formula in Lemma \ref{hfixpt} reads
\begin{eqnarray*}
h(G/P,\mtr {\cal L}_\lambda)&=&\sum_{w\in
W_G/W_K}\frac{1}{\prod_{\a\in\Psi}w\a(Y)}
\\&&\cdot
\sum_{l=1}^{n+1}
\sum_{\a\in\Psi}
\frac{(w\lambda(Y))^{n+1}-(w\lambda(Y))^{n+1-l}({S_{w\a}
w\lambda(Y)})^l}{2lw\a(Y)}\,\,.
\end{eqnarray*}
A good choice for $Y$ is the dual of $\rho$.
\beginProof 
Applying the Bott residue formula to formula (\ref{harmo}) yields
$$
h(G/P,\mtr {\cal L}_\lambda)=\sum_{p\in M^Y}\sum_{l=0}^n{n+1\choose l+1}
\frac{\phi^{n-l}}{2(l+1)\prod\theta_\nu}
\sum_{j\in\BN}\sum_{\theta_\nu\in\Psi_j}
(-\theta_\nu)^l j^{l+1}
$$
(alternatively, one could apply \cite[Th. 11]{K2} to $f_j(m,k)$ from equation
(\ref{hupp}) and proceed as in  \cite[proposition 3.7]{KR2}).
Using the formula $$\sum_{l=1}^n {n\choose l}\frac{x^l}l= \sum_{l=1}^n
\frac{(1+x)^l-1}l\qquad (x\in\BR)\,\,,$$ we get the desired result.
\endProof
{\bf Example. }We shall express the height of the Grassmannian $G(m,k)$ with
$G(m,k)(\BC)={\bf U}(m)/{\bf U}(k)\times {\bf U}(m-k)$ using Lemma
\ref{hfixpt}. In this case, the Weyl group of
$G$ is the permutation group $S_m$ of $m$ elements and the fixed point set can be
identified with $S_m/S_k\times S_{m-k}$. Let $(\epsilon_\nu)_{\nu=1}^m$ be the
cartesian base of $\BC^m$, which we identify with ${\frak t}^\vee$ as in
\cite[section V.6]{BtD}. A short look to the classical tables of roots
(e.g. \cite[VI, planche I]{Bour} or \cite[Proposition V.6.2]{BtD}) reveals that
$\Psi=\{\epsilon_\mu-\epsilon_\nu|1\leq\mu\leq k<\nu\leq m\}$.
 Also, there is a unique positive primitive line bundle
${\cal L}_\lambda$ on
$G(m,k)$ with $\lambda=\sum_{\nu=1}^k \epsilon_k$. In particular, $\Psi=\Psi_1$ . Set
$I_m:=\{1,\dots,m\}$; when considering the action corresponding to $Y=\sum_{\nu=1}^m
\nu \epsilon_\nu^*$, we get by Lemma
\ref{hfixpt}
$$
h(G(m,k),\mtr {\cal L}_\lambda)=\sum_{I\subseteq I_m\atop \#I=k}
\frac{(\sum_{c\in I}c)^{k(m-k)+1}}{\prod_{a\in I\atop b\in I_m\setminus I}(a-b)}
\sum_{a\in I\atop b\in I_m\setminus I} \sum_{l=1}^{k(m-k)+1}
\frac{1-(1-\frac{a-b}{\sum_{c\in I}c})^l}{2 l (a-b)}\,\,.
$$

{\bf Example. }Assume that $G/P$ is embedded as a hypersurface of degree $d$ into
${\bf P}^{n+1}$ via ${\cal L}_\lambda$. 
Let
$N$ be its normal bundle. Classically, there is an exact sequence
$$
0\to{\cal O}\to\bigoplus_1^{n+2}{\cal O}(1)\to T{\bf P}^{n+1}\to 0
$$
\cite[Ex. 8.20.1]{Hart}. Furthermore, by the adjunction formula (\cite[p. 146]{GH})
$$
N={\cal O}(d)_{|M}\,\,.
$$
Thus $\ch(TM)=\ch(T{\bf P}^{n+1}_{|M})-\ch(N)=(n+2)e^{c_1({\cal O}(1))}-1 -e^{d
c_1({\cal O}(1))}$. Assume now that $E_1=TM$. Then equation (\ref{harmo}) yields
\begin{eqnarray}
h(G/P,\mtr {{\cal O}(1)})&=&
\sum_{l=1}^n \frac{(-1)^l}{2(l+1)}{n+1\choose l+1}\int_M (n+2-d^l) c_1({\cal
O}(1))^n
\nonumber\\&&
+\frac{n+1}2\int_M n c_1({\cal O}(1))^n
\nonumber\\&=&
\sum_{l=1}^n \frac{d (-1)^l}{2(l+1)}{n+1\choose l+1}(n+2-d^l)+\frac{n(n+1)d}2
\nonumber\\&=&\label{topo}
\sum_{l=2}^{n+1} \frac1{2l}\left[d(n+2)-1+(1-d)^l\right]
\end{eqnarray}
By \cite[Th. III.2.3 and next paragraph]{Ko}, all complex hypersurfaces
of dimension $n>1$ with a non-trivial holomorphic $\BC^*$-action have degree $d\leq2$.
We reobtain the formula for the height of ${\bf P}^n$ \cite{BoGS}
$$
h({\bf P}^n,\mtr {{\cal O}(1)})=\frac{n+1}2\sum_{k=1}^n\frac1{k}-\frac n{2}\,\,.
$$
Also \cite[VI, planche IV]{Bour} shows $\Psi=\Psi_1$ for the
even-dimensional smooth quadric $Q_{2m}$ with $Q_{2m}(\BC)={\bf SO}(2m+2)/{\bf
SO}(2m){\bf SO}(2)$. Thus  we reobtain the result from
\cite[Cor. 2.2.10]{CM}
$$
h(Q_{2m},\mtr {{\cal O}(1)})=(2m+1)\sum_{k=1}^{2m-1}\frac1{k}+\frac1{2}
\sum_{k=1}^{m-1}\frac1{k}-2m+1+\frac1{m}\,\,.
$$
Now consider the case of the odd-dimensional smooth quadric $Q_{2m-1}$ with
$Q_{2m-1}(\BC)={\bf SO}(2m+1)/{\bf SO}(2m-1){\bf SO}(2)$. Identifying $\BC^m$ with
$\frak t^\vee$ as in \cite[proposition V.6.5]{BtD} 
and denoting the cartesian base by
$(\epsilon_\nu)_{\nu=1}^m$, we notice $\Psi=\{\epsilon_1\pm \epsilon_\nu|1<\nu\leq
m\}\cup\{\epsilon_1\}$ and $\lambda=\epsilon_1$ (compare also \cite[VI,planche
II]{Bour}). Thus $TM=E_1+E_2$ with $E_2={\cal O}(1)$.
Hence there is an additional term to equation (\ref{topo}) given by
\begin{eqnarray*}
\lefteqn{
\sum_{l=0}^n \frac{(-1)^l}{2(l+1)}{n+1\choose l+1}\int_M (2^{l+1}-1) c_1({\cal
O}(1))^n
}\\&=&
\sum_{l=0}^n \frac{d(-1)^l}{2(l+1)}{n+1\choose l+1} (2^{l+1}-1)
=
\sum_{l=1}^{n+1}\frac{(-1)^{l+1}}l
\end{eqnarray*}
and we find
$$
h(Q_{2m-1},\mtr {{\cal O}(1)})=(2m+1)\sum_{k=1}^{2m-1}\frac1{k}-\frac1{2}
\sum_{k=1}^{m-1}\frac1{k}-2m+1
$$
which is exactly the same value as obtained in \cite[Cor. 2.2.10]{CM} for the singular
quadric $X_0^2+\cdots+X^2_{2m}=0$.

\end{document}